\documentstyle[12pt]{article}
\newlength{\mytopmargin}
\newlength{\myleftmargin}
\begin{document}
\vspace{3.5cm}
\noindent
\begin{center}{  \bf 
GENERALIZED WEIGHT FUNCTIONS AND THE \\
 MACDONALD POLYNOMIALS} 
\end{center}
\vspace{5mm}

\noindent
{\center T.H.~Baker\footnote{email: tbaker@maths.mu.oz.au; supported by the ARC}
and P.J.~Forrester\footnote{email: matpjf@maths.mu.oz.au; supported by the ARC}\\
\it Department of Mathematics, University of Melbourne, Parkville, Victoria
3052, Australia}
\vspace{.5cm}

\small
\begin{quote}
A weight function which $q$-generalizes the  ground state wave function 
of the multi-component Calogero-Sutherland
quantum many body system is introduced. Conjectures, and some proofs in
special cases, are given for a constant term identity involving this function.
A Gram-Schmidt procedure with respect to the inner product associated with
 the weight function is used to define
orthogonal polynomials in one of the components, which are 
conjectured to be the Macdonald polynomials 
$P_\kappa(w_1,\dots,w_{N_0};qt^p,t)$,
and a proof is given in a special case. Conjectures are also given for an
adjoint property of the Macdonald operator with respect to the inner product
associated with the weight function, and the normalization of the Macdonald
polynomial with respect to the same inner product.
\end{quote}
\vspace{.5cm}
\noindent
{\bf 1. INTRODUCTION}

\vspace{3mm}
In two recent studies \cite{forr1,forr2} of the multi-component
Calogero-Sutherland model (quantum many body system with $1/r^2$ pair
potential), one of us has been led  to formulate a number of conjectures
concerning Jack polynomials \cite{stan,macd} and the function
\renewcommand{\theequation}{1.1}
\begin{eqnarray}
\lefteqn{|\psi_0(\{z_j^{(\alpha)}
\}_{\alpha = 1,\dots,p \atop j=1,\dots,N_\alpha},
\{w_j\}_{j=1,\dots,N_0})|^2}
\nonumber \\
& := & \prod_{\alpha=1}^p \prod_{1 \le j < k \le N_{\alpha}}
|z_j^{(\alpha)} - z_k^{(\alpha)}|^{2\lambda + 2}
\prod_{1 \le j' < k' \le N_0}|w_{k'}-w_{j'}|^{2\lambda}\nonumber \\
& &\times \prod_{1 \le \alpha < \beta \le p}
\prod_{j=1}^{N_\alpha}\prod_{k=1}^{N_\beta}
|z_j^{(\alpha)} - z_k^{(\beta)}|^{2 \lambda }\:
\prod_{\alpha=1}^p\prod_{j=1}^{N_\alpha}\prod_{j'=1}^{N_0}
|z_j^{(\alpha)} - w_{j'}|^{2 \lambda },
\end{eqnarray}
where $w_j := e^{2 \pi i y_j}$ and $z_j^{(\alpha)} :=
 e^{2 \pi i x_j^{(\alpha)}}$, which is the absolute value squared of the ground
state wave function. In fact it appears that the Jack polynomials can be
constructed via a Gram-Schmidt procedure based on (1.1) as a weight function.

Explicitly, define an inner product by
$$
\langle f|g \rangle_{N_0,\dots,N_p;\lambda} := 
  \prod_{l=1}^{N_0} \int_{-1/2}^{1/2} d y_l \,
\prod_{\alpha = 1}^p \prod_{l=1}^{N_\alpha}
 \int_{-1/2}^{1/2} d x_l^{(\alpha)} \,|\psi_0(\{z_j^{(\alpha)}
\}_{\alpha = 1,\dots,p \atop j=1,\dots,N_\alpha},
\{w_j\}_{j=1,\dots,N_0})|^2 f^* \, g.
\eqno (1.2)
$$
(this notation differs from that used in ref.~\cite{forr2} in that the weight
function is explicitly included in the r.h.s.). Let $\kappa$ denote a partition
and define  a symmetric polynomial
in the variables $w_1, \dots, w_{N_0}$, denoted $p_\kappa(w_1, \dots, w_{N_0})$,
by the following properties:

(i) $p_\kappa(w_1, \dots, w_{N_0}) = m_\kappa + \sum_{\mu < \kappa} a_\mu m_\mu$,
where $|\mu| = |\kappa|$, $\mu < \kappa$ is with respect to reverse lexicographical ordering of
the partitions, $m_\mu$ refers to the monomial symmetric function with
exponents $\mu = (\mu_1, \dots, \mu_N) $ in the variables
$w_1, \dots, w_{N_0}$ and $a_\mu$ is the corresponding coefficient;

(ii) for all $N_1, \dots, N_p \ge \kappa_1 - 1$, $\langle  p_\kappa
|  p_\sigma \rangle_{N_0,\dots,N_p;\lambda} = 0$ for $\kappa \ne \sigma$.

\vspace{.2cm}
\noindent
(in ref.~\cite{forr2} condition (ii) required $N_1, \dots, N_p \ge \kappa_1$;
this was weakened to the above statement in Conjecture 2.2 of
ref.~\cite{forr2}).
 Then, according to Conjecture 2.4 of ref.~\cite{forr2}, the
polynomials
$p_\kappa$ are given in terms of the Jack polynomials by
$$
p_\kappa(w_1, \dots, w_{N_0}) =  J_\kappa^{(1 + 1/\lambda)}
(w_1, \dots,w_{N_0})
\eqno (1.3)
$$
(here the normalization of $J_\kappa$ is chosen so that the coefficient of 
$m_\kappa$ is unity).

Conjecture 2.4 of ref.~\cite{forr2} is to be contrasted with the known theorem
\cite{stan,macd} for the 
construction of the Jack polynomials via a Gram-Scmidt procedure based on the
$p=0$ cases of the weight function (1.1) and the inner product (1.2). Then the
symmetric polynomials with properties (i) and (ii) above (in (ii), since $p=0$,
there is no restriction on 
$\kappa_1$) are given by 
$$
p_\kappa(w_1, \dots, w_{N_0}) =  J_\kappa^{(1/\lambda)}
(w_1, \dots,w_{N_0}).
\eqno (1.4)
$$
Note that the parameter of the Jack polynomial here is $1/\lambda$.

The theory of Jack polynomials has been $q$-generalized by Macdonald [4] to
give a theory of what are now referred to as Macdonald polynomials. This has
motivated us to seek $q$-generalizations of the conjectures (and some theorems)
contained in refs.~\cite{forr1,forr2}. We begin in Section 2 by $q$-generalizing
the weight function (1.1) and considering the $q$-generalizations of conjectured
constant term identities given in ref.~\cite{forr1}.
In Section 3 proofs of the conjectures
of Section 2 are provided in certain cases. The Gram-Schmidt procedure is used
with the
$q$-generalization of (1.1) to define $q$-generalizations of the polynomials
(1.3) in Section 4, and a conjecture  is given relating these
polynomials to the Macdonald polynomials. In fact we are led to conclude
that (1.3) is only correct for $p=1$. For general $p$ we have new evidence
which suggests that Conjecture 2.4 of ref.~\cite{forr2} should read
$$
p_\kappa(w_1, \dots, w_{N_0}) =  J_\kappa^{(p+1/\lambda)}
(w_1, \dots,w_{N_0}).
\eqno (1.5)
$$
We also provide a conjecture for a
normalization integral involving the Macdonald polynomials, generalizing the
conjecture given in ref.~\cite{forr2} in the Jack polynomial case. In Section 5
some proofs of the conjectures of Section 4 are provided in certain cases, while
our results are briefly summarized in Section 6. In the Appendix we use a known
generalization of the so called $q$-Morris theorem (see e.g. ref.~\cite{zeil}),
which is the
$p=0$ case of the constant term identities considered in Section 2, to derive
the expansion of the power sums in terms of Macdonald polynomials.

\vspace{.5cm}
 \noindent
{\bf 2. THE $q$-GENERALIZED WEIGHT FUNCTION AND CONSTANT TERM IDENTITIES}

\vspace{.3cm}
\noindent
{\bf 2.1 Revision of the case $p=0$}

\noindent
For $p=0$ and $\lambda$ integer, since $|w_j| = 1$, (1.1) can be written as
$$
|\psi_0(
\{w_j\}_{j=1,\dots,N_0})|^2 = 
\prod_{1 \le j < k \le N_0} \Big (1 - {w_k \over w_j}\Big )^\lambda
 \Big (1 - {w_j \over w_k}\Big )^\lambda.
\eqno (2.1)
$$
This was first $q$-generalized by Andrews \cite{andr} as
$$
|\psi_0(
\{w_j\}_{j=1,\dots,N_0};q)|^2 := 
\prod_{1 \le j < k \le N_0} \Big (q{w_k \over w_j};q\Big )_\lambda
 \Big ({w_j \over w_k};q\Big )_\lambda
\eqno (2.2)
$$
where
$$
(a;q)_\lambda := \prod_{l=0}^{\lambda - 1} (1 - a q^l), \quad \lambda \in
Z_{\ge 0}.
\eqno (2.3)
$$

The criterium used to choose this $q$-generalization (note that unlike (2.1),
(2.2) is not symmetric in $w_1, \dots, w_{N_0}$), additional to requiring that
(2.2) reduces to (2.1) when $q=1$, was that the Dyson identity \cite{dyso}
$$
{\rm CT}\, |\psi_0(
\{w_j\}_{j=1,\dots,N_0})|^2 = {(\lambda N_0)! \over \lambda !^N_0}, \quad
\lambda \in Z_{\ge 0},
\eqno (2.4)
$$
generalizes as
$$
{\rm CT}\, |\psi_0(
\{w_j\}_{j=1,\dots,N_0};q)|^2 = {\Gamma_q (\lambda N_0 + 1) \over
(\Gamma_q(\lambda + 1))^N_0}, \quad \lambda \in Z_{\ge 0}
\eqno (2.5)
$$
where
$$
\Gamma_q(n + 1) := \prod_{j=1}^n {1 - q^j \over 1 - q}
\eqno (2.6)
$$
and CT denotes the constant term in the Laurent polynomial.
Note the restriction $\lambda \in Z_{\ge 0}$ in the above formulas. For general
$\lambda$ we interpret (2.3) as
$$
(a;q)_\lambda := {(a;q)_\infty \over (aq^\lambda;q)_\infty}.
\eqno (2.7)
$$
As pointed out by Stembridge \cite{stem}, the identity (2.5) still holds with
$$
\Gamma_q(x) := {(q;q)_\infty \over (q^{x};q)_\infty} (1-q)^{1-x}
= {(q;q)_{x-1} \over (1 - q)^{x-1}}
\eqno (2.8)
$$
and
$$
{\rm CT} \, f(w_1, \dots, w_N) = \prod_{l=1}^N \int_{-1/2}^{1/2} dx_l \,
f(e^{2 \pi i x_1}, \dots, e^{2 \pi i x_N}).
\eqno (2.9)
$$
\vspace{.2cm}
\noindent
{\bf 2.2 $q$-generalization for general $p$}

\noindent
Motivated by the $q$-generalization (2.2) of (2.1), we formulated the 
$q$-generalization of (1.1) as
\renewcommand{\theequation}{2.10}
\begin{eqnarray}
\lefteqn{|\psi_0(\{z_j^{(\alpha)}
\}_{\alpha = 1,\dots,p \atop j=1,\dots,N_\alpha},
\{w_j\}_{j=1,\dots,N_0};q)|^2}
\nonumber \\
& = & \prod_{\alpha=1}^p \prod_{1 \le j < k \le N_{\alpha}}
\Big ({z_j^{(\alpha)}\over z_k^{(\alpha)}};q \Big )_{\lambda + 1}
\Big (q{z_k^{(\alpha)}\over z_j^{(\alpha)}};q \Big )_{\lambda + 1}
\prod_{1 \le j' < k' \le N_0}\Big ({w_{j'}\over w_{k'}};q
\Big )_{\lambda}
\Big (q{w_{k'}\over w_{j'}};q
\Big )_{\lambda}
\nonumber
\\ & &\times \prod_{1 \le \alpha < \beta \le p}
\prod_{j=1}^{N_\alpha}\prod_{k=1}^{N_\beta}
\Big ({z_j^{(\alpha)} \over z_k^{(\beta)}};q \Big )_{ \lambda }
\Big (q{z_k^{(\beta)} \over z_j^{(\alpha)}};q \Big )_{ \lambda }
\:
\prod_{\alpha=1}^p\prod_{j=1}^{N_\alpha}\prod_{j'=1}^{N_0}
\Big ({z_j^{(\alpha)} \over w_{j'}};q \Big )_{ \lambda }
\Big (q{w_{j'} \over z_{j}^{(\alpha)}};q \Big )_{ \lambda }.
\end{eqnarray}
Indeed this $q$-generalization appears to generalize an integration formula for
(1.1), conjectured in ref.~\cite{forr1}, in the same way that (2.5) generalizes
(2.4).

To be more explicit, let us consider the case $p=1$, and introduce the notation
\renewcommand{\theequation}{2.11}
\begin{eqnarray}\lefteqn{D_p(N_1;N_0;a,b,\lambda)}
\nonumber \\&:= & 
\left ( \prod_{l=1}^{N_0} \int_{-1/2}^{1/2} d x_l \, w_l^{ (a-b)/2}
|1+w_l|^{a+b} \right )
\left ( \prod_{l=1}^{N_1}
 \int_{-1/2}^{1/2} d x_l \, z_l^{ (a-b)/2}
|1+z_l|^{a+b} \right )
 \nonumber \\
& & \times|\psi_0(\{z_j
\}_{ j=1,\dots,N_1},
\{w_j\}_{j=1,\dots,N_0})|^2. 
\end{eqnarray}
In ref.~\cite[eq.(3.21)]{forr1} it was conjectured that
\renewcommand{\theequation}{2.12}
\begin{eqnarray}
\lefteqn{D_1(N_1;N_0;a,b,\lambda)} \nonumber \\ & = &
\prod_{j=0}^{N_1 - 1} {(j+1) \Gamma ((\lambda + 1)j+a+b+\lambda N_0+1)
\Gamma((\lambda + 1)(j+1)+\lambda N_0) \over
\Gamma(1 + \lambda)  \Gamma ((\lambda + 1)j+a+\lambda N_0+1)
\Gamma ((\lambda + 1)j+b+\lambda N_0+1)}\nonumber \\
& & \times \prod_{l=0}^{N_0-1} {\Gamma (a+b+1+\lambda l) \Gamma (1+\lambda(l+1))
\over
\Gamma (a+1+\lambda l)\Gamma (b+1+\lambda l) \Gamma (1+\lambda)}
\end{eqnarray}
To formulate the $q$-generalization of (2.11) we note that for $a$ and $b$
integers and $|u| = 1$
$$
u^{(a-b)/2} |1 + u|^{a+b} = (1+u)^a(1+1/u)^b.
\eqno (2.13)
$$
This suggests we define
\renewcommand{\theequation}{2.14}
\begin{eqnarray}\lefteqn{D_1(N_1;N_0;a,b,\lambda;q)}
\nonumber \\ 
&:= & \prod_{l=1}^{N_0} \int_{-1/2}^{1/2} d y_l \, 
(-w_l;q)_a \Big (-{q \over w_l};q \Big )_a 
 \prod_{l=1}^{N_1}
 \int_{-1/2}^{1/2} d x_l \, 
(-z_l;q)_b  \Big (-{q \over z_l};q \Big )_b 
 \nonumber \\
& & \times|\psi_0(\{z_j
\}_{ j=1,\dots,N_1},
\{w_j\}_{j=1,\dots,N_0};q)|^2  \nonumber \\
& = & \prod_{l=1}^{N_0} \int_{0}^{1} d y_l \,
(w_l;q)_a \Big ({q \over w_l};q \Big )_a 
 \prod_{l=1}^{N_1}
 \int_{0}^{1} d x_l \, 
(z_l;q)_b  \Big ({q \over z_l};q \Big )_b 
 \nonumber \\
& & \times \prod_{1 \le j < k \le N_1}
\Big ({z_j \over z_k};q \Big )_{\lambda + 1}
\Big (q{z_k \over z_j};q \Big )_{\lambda + 1}
\prod_{1 \le j' < k' \le N_0}\Big ({w_{j'}\over w_{k'}};q
\Big )_{\lambda}
\Big (q{w_{k'}\over w_{j'}};q
\Big )_{\lambda} \nonumber \\
& & \times \prod_{j = 1}^{N_1} \prod_{j' = 1}^{N_0} \Big ({z_j \over
w_{j'}};q
\Big )_{
\lambda }
\Big (q{w_{j'} \over z_{j}};q \Big )_{ \lambda }
\end{eqnarray}
On the basis of exact computer generated data  (obtained for $\lambda = 1$ and
2, with various `small' values of $N_0$, $N_1$, $a$ and $b$) and some analytic
evaluations in certain special cases (presented in the next section) we
make the following conjecture for the evaluation of (2.11).

\vspace{.2cm}
\noindent
{\bf Conjecture 2.1} \quad We have
\renewcommand{\theequation}{2.15}
\begin{eqnarray}
\lefteqn{D_1(N_1;N_0;a,b,\lambda)} \nonumber \\ & = &
{\Gamma_{q^{\lambda + 1}}(N_1 + 1) \over (\Gamma_q(1 + \lambda))^{N_0 + N_1}}
\prod_{j=0}^{N_1 - 1} { \Gamma_q ((\lambda + 1)j+a+b+\lambda N_0+1)
\Gamma_q((\lambda + 1)(j+1)+\lambda N_0) \over
  \Gamma_q ((\lambda + 1)j+a+\lambda N_0+1)
\Gamma_q ((\lambda + 1)j+b+\lambda N_0+1)}\nonumber \\
& & \times \prod_{l=0}^{N_0-1} {\Gamma_q (a+b+1+\lambda l) \Gamma_q
(1+\lambda(l+1))
\over
\Gamma_q (a+1+\lambda l)\Gamma_q (b+1+\lambda l) }
\end{eqnarray}
Note that the `base' of the $q$-gamma function in the denominator of the first
term is $q^{\lambda + 1}$ whereas in all other terms it is $q$. Also, when
$N_1 = 0$, note that this reduces to the so-called $q$-Morris theorem (see 
e.g.~ref.~\cite{zeil}).

\vspace{.2cm}
\noindent
{\bf 2.3 The $q$-generalized integral for general $p$}

\noindent
In the $q=1$ case it was conjectured \cite[eq.(4.8a)]{forr1} that for general $p$
the analogue of the integral (2.14), $D_p$ say, satisfies a functional equation.
Using this functional equation it  was shown $D_p$ can be uniquely determined
by a recurrence. To $q$-generalize this result, let
\renewcommand{\theequation}{2.16a}
\begin{eqnarray}\lefteqn{D_p(N_1, \dots,N_p;N_0;a,b,\lambda;q)}
\nonumber \\& := & 
 \prod_{l=1}^{N_0} \int_{-1/2}^{1/2} d y_l  \prod_{\alpha =
1}^p
\prod_{l=1}^{N_\alpha}
 \int_{-1/2}^{1/2} d x_l^{(\alpha)}  \,
A(\{- z_j^{(\alpha)} \}_{\alpha = 1,\dots,p \atop j=1,\dots,N_\alpha},
 \{- w_j\}_{j=1,\dots,N_0};q)
 \nonumber \\
& & \times|\psi_0(\{z_j^{(\alpha)}
\}_{\alpha = 1,\dots,p \atop j=1,\dots,N_\alpha},
\{w_j\}_{j=1,\dots,N_0};q)|^2 
\end{eqnarray}
where
\renewcommand{\theequation}{2.16b}
\begin{equation}
A(\{ z_j^{(\alpha)}\}_{\alpha = 1,\dots,p \atop j=1,\dots,N_\alpha},
 \{ w_j\}_{j=1,\dots,N_0};q) :=
\prod_{l=1}^{N_0} (w_l;q)_a \Big ({q \over w_l};q \Big )_a
\prod_{\alpha = 1}^p \prod_{j=1}^{N_\alpha} 
(z_l^{(\alpha)};q)_b  \Big ({q \over z_l^{(\alpha)}};q \Big )_b.
\end{equation}

Guided by the conjecture in ref.~\cite[eq.(4.8a)]{forr1},
and Conjecture 2.1 above, we can
make the following conjecture for the general $p$ case.

\vspace{.2cm}
\noindent
{\bf Conjecture 2.2} \quad
For $N_{p} \ge N_j - 1$ $(j = 1, \dots, p-1)$ we have
\renewcommand{\theequation}{2.17}
\begin{eqnarray}
\lefteqn{{D_p(N_1, \dots,N_{p-1},N_p + 1;N_0;a,b,\lambda;q) \over
D_p(N_1, \dots,N_{p-1},N_p;N_0;a,b,\lambda;q)} = {[N_p+1]_{q^{\lambda +1}}
\over \Gamma_q(\lambda + 1)}} \nonumber \\
& &  \times { \Gamma_q ((\lambda + 1)N_p+a+b+\lambda \sum_{j=0}^{p-1}N_j+1)
\Gamma_q((\lambda + 1)(N_p+1)+\lambda \sum_{j=0}^{p-1}N_j ) \over
  \Gamma_q ((\lambda + 1)N_p+a+\lambda \sum_{j=0}^{p-1}N_j+1)
\Gamma_q ((\lambda + 1)N_p+b+\lambda \sum_{j=0}^{p-1}N_j+1)}.
\end{eqnarray}

\vspace{.2cm}
\noindent
In Conjecture 2.2 we have introduced the notation
$$
[a]_q := {1 - q^a \over 1 - q}.
\eqno (2.18)
$$
In the limit $q \rightarrow 1$ the formula in Conjecture 2.2 is equivalent
to the functional equation conjectured in ref.~\cite[eq.(4.8a)]{forr1}.

\vspace{.5cm}
 \noindent
{\bf 3. ANALYTIC EVALUATIONS OF THE $q$-GENERALIZED INTEGRAL}

\vspace{.3cm}
\noindent
In this section Conjecture 2.1 will be proved in some special cases.

\vspace{.2cm}
\noindent
{\bf 3.1 The case $p=1$, $a = \lambda$ and general $N_0,N_1,b,\lambda$}

\noindent
It was noted in ref.\cite{forr1} that the integral evaluation in Conjecture 2.1
for
$q=1$ in the case $p=1$, $a=b=0$ and general $N_0,N_1,\lambda$ follows from a
theorem of Bressoud and Goulden \cite{brgo2}. This theorem has a $q$-counterpart,
obtained by the same authors in an earlier publication.

\vspace{.2cm}
\noindent
{\bf Proposition 3.1} \cite[Proposition 2.4, with $A$ replaced by
$\bar{A}$, the complement of $A$, to be consistent with the formulation in
ref.~\cite{brgo1}]{brgo2}
\quad Let $a_1, \dots, a_N$ be positive integers, $A$ be an arbitrary subset of
$\{ (i,j): 1 \le i < j \le n \}$, $G_A$ be the set of permutations $\sigma$ on
$\{1, \dots, n\}$ (with $\sigma(i) := \sigma_i$) whose inversions 
$
I(\sigma) := \{ (\sigma_i, \sigma_j) : j > i \quad {\rm and}
\quad \sigma_j < \sigma_i \}
$
are contained in $A$:
$$
G_A = \{ \sigma: \:{\rm if} \: (j,i) \in I(\sigma), \: {\rm then} \: (i,j)
\in A \},
$$
and let $\chi (T)$ be the characteristic function which is 1 if $T$ is true,
0 otherwise. We have
$$
{\rm CT} \, \prod_{1 \le i < j \le n} \Big ( q{x_i \over x_j};q
\Big )_{a_i} \Big ( {x_j \over x_i};q \Big )_{a_j - \chi((i,j) \notin A)}
= {\Gamma_q(a_1 + \dots + a_n + 1) \over
\Gamma_q(a_1) \dots \Gamma_q(a_n)}
S_n(\{a_j\}_{j=1, \dots, n};G_A),
$$
where
$$
S_n(\{a_j\}_{j=1, \dots, n};G_A) :=
\sum_{\sigma \in G_A} q^{\sum_{(i,j) \in I(\sigma)} a_i} 
\prod_{l=1}^n { 1 - q \over 1 - q^{a_{\sigma_1} + \dots + a_{\sigma_l}}}.
$$

\vspace{.2cm}
In this theorem, suppose $n = 1+ N_0 + N_1$,
$$
A = \{(i,j): 1 \le i < j \le N_0 + 1 \quad {\rm or} \quad
N_1 + 2 \le i < j \le N_0 + N_1 +1 \},
$$
$$
a_1 = b, \quad, a_2 = \dots = a_{N_0 + 1} = \lambda, \quad
a_{N_0 +2} = \dots a_{N_0 + N_1 + 1} = \lambda + 1,
\eqno (3.1)
$$
and replace $x_2, \dots,x_{N_0 + N_1 +1}$ by $x_1x_2, \dots, x_1x_{N_0 + N_1
+1}$. We see that the l.h.s. is of the form (2.14) with $a = \lambda$ and
$b, \lambda$ arbitrary positive integers. Proposition 3.1 therefore gives
$$
 D_1(N_1;N_0;\lambda,b,\lambda;q) = {\Gamma_q((\lambda + 1)N_1 + \lambda N_0 +
b+1) 
\over \Gamma_q(b ) (\Gamma_q(\lambda))^{N_0}(\Gamma_q(\lambda +
1))^{N_1}} S_{N_0,N_1}(\{a_j\}_{j=1, \dots, n};G_A).
\eqno (3.2)
$$
Our proof of Conjecture 2.1 for the evaluation of $
D_1(N_1;N_0;\lambda,b,\lambda;q)$ now follows from the following evaluation of
$S_{N_0,N_1}$.

\vspace{.2cm}
\noindent
{\bf Proposition 3.2} \quad With $A$, $\{a_j\}$ given by (3.1) we have
$$
 S_{N_0,N_1}(\{a_j\};G_A) = 
{1 \over [b]_q [\lambda]_q^{N_0}[\lambda + 1]_q^{N_1}}
\prod_{j=1}^{N_1} {[(\lambda + 1)j]_q \over
[(\lambda + 1)j + \lambda N_0 + b]_q }.
$$

\vspace{.2cm}
\noindent
{\bf Proof} \quad We will establish a recurrence relation in $N_1$. To solve the
recurrence we require the value of $S_{N_0,0}$. This is obtained by noting that
when $N_1=0$, $D_1$ is given by the $q$-Morris theorem (Conjecture 2.1 with
$N_1 = 0$). Comparison with (3.2) then  gives
$$
S_{N_0,0}(\{a_j\};G_A) = 
{1 \over [b]_q [\lambda]_q^{N_0}}.
\eqno (3.3)
$$
The recurrence is obtained by noting that all permutations in $G_A$ are of the
form $\sigma = (\sigma', \sigma'')$, where $\sigma'$ is a permutation of
$\{1, \dots, N_0 + 1\}$, and $\sigma''$ is a permutation of $\{N_0 + 2, \dots,
N_0 + N_1 + 1\}$. Thus, if $N_1$ is increased by 1 only $\sigma''$ can be
affected. Furthemore, for each $\sigma = (\sigma', \sigma'')$ in $G_A$ before
increasing $N_1$ by 1, there are $N_1 + 1$ permutations in $G_A$ after $N_1$ is
increased by 1, which are given by $\sigma  = (\sigma', \sigma''|_{k \mapsto
N_0 + N_1 + 2},k)$, $k = N_0 + 2, \dots, N_0 + N_1 + 2$ (for $k = N_0 + N_1 +
2$, $\sigma''$ remains unchanged). Denote these permutations by $G_{A}(k)$ so
that for $N_1$ increased by 1, $G_A = \cup_{k=N_0 + 2}^{N_1 + N_0 + 2}
G_{A}(k)$.  The facts that the replacement $k \mapsto N_0 + N_1 + 2$ in
$\sigma''$ creates $N_0 + N_1 + 2 - k$ new inversions, that $a_i = \lambda + 1
$ for $i = N_0 + 2, \dots, N_0 + N_1 + 2$, and that
$$
(1-q)/(1-q^{a_{\sigma_1 }+ \dots + a_{\sigma_{N_0 + N_1 + 2}}}) =
(1-q)/(1-q^{a_1 + \dots + a_{N_0 + N_1 + 2}})
$$
 is a common factor in the summand now gives the recurrence
$$
 S_{N_0,N_1+1}(\{a_j\};G_A(k)) = 
{q^{(\lambda + 1)(N_0 + N_1 + 2 - k)}\over
 [(\lambda + 1)(N_1 + 1) + \lambda N_0
+ b]_q} S_{N_0,N_1}(\{a_j\};G_A(k))
$$
Summing over $k$ we have
$$
 S_{N_0,N_1+1}(\{a_j\};G_A) = 
{[(\lambda + 1)(N_1 + 1)]_q \over [(\lambda + 1)(N_1 + 1) +
\lambda N_0 + b]_q [\lambda + 1]_q}
 S_{N_0,N_1}(\{a_j\};G_A),
$$
which upon iteration and use of the initial condition gives the stated result.

\vspace{.2cm}

Substituting the result of Proposition 3.2 in (3.2) evaluates
$D_1(N_1;N_0;\lambda,b,\lambda;q)$ for $\lambda$ and $b$ arbitrary positive
integers; by a simple lemma of Stembridge \cite[lemma 3.2]{stem}, the validity of
the positive integer case implies the validity for all (complex) $\lambda$ and
$b$. Comparison between the resulting expression for $D_1$ and the expression
of Conjecture 3.1 in the case $a=\lambda$ shows, after simplification of the
latter, that the two expressions are identical.

\vspace{.2cm}
\noindent
{\bf 3.2 The case $p=1$, $N_1 = 2$ and general $N_0,N_1,a,b,\lambda$}

\setcounter{equation}{3}
\setcounter{section}{3}
\renewcommand{\theequation}{\thesection.\arabic{equation}}
\newcommand{\la}{\lambda}
\newcommand{\si}{\sigma}
\newcommand{\al}{\alpha}
\newcommand{\Ga}{\Gamma}
\newcommand{\ga}{\gamma}
\newcommand{\ka}{\kappa}
\newcommand{\de}{\delta}
\def\rlx{\relax\leavevmode}
\def\IZ{\rlx\hbox{\small \sf Z\kern-.4em Z}}
\def\IR{\rlx\hbox{\rm I\kern-.18em R}}
\def\IN{\rlx\hbox{\rm I\kern-.18em N}}
\def\ID{\rlx\hbox{\rm I\kern-.18em D}}
\def\IC{\rlx\hbox{\,$\inbar\kern-.3em{\rm C}$}}

We address the case $N_1=2$, as this is the first non-trivial case;
when $N_1=1$, it is clear that $D_1(1;N_0;a,b,\la;q)=D_0(N_0+1;a,b,\la;q)$.
To prove Conjecture 2.1 in this particular case, we adopt the
method of Stembridge and Stanton\cite{stem}. The essence of this method,
applied to the problem at hand, is to express the constant term 
of the two-component function in terms of the constant term of
the one-component function (that is, the function appearing in the
$q$-Morris identity), by means of a partial expansion of the 
two-component function. 
\bigskip

It will prove useful to briefly summarize the results of
Zeilberger \cite{zeil} as we shall be aiming to extend his proof
given there, of the $q$-Morris identity. Let
$[x^0]f(x_1,\ldots,x_n):= {\rm CT} f(x_1,\ldots,x_n)$, and in
general let $[x^\beta]f(x_1,\ldots,x_n)$ denote the
coefficient of $x^{\beta}:=x_1^{\beta_1}x_2^{\beta_2}\cdots x_n^{\beta_n}$
in the expansion of $f$. For notational convenience, this can be
extended so that for a general function $g(x)=\sum_{\beta}a_{\beta}
x^{\beta}$, one writes $[g]f := \sum_{\beta}a_{\beta}[x^{\beta}]f$.

The ``reduced'' $q$-Morris identity takes the form
\begin{eqnarray}
[x^0]F_0(x) &=& \prod_{j=0}^{n-1}\frac{\Ga_q(\la j+a+b+1)\Ga_q(\la(j+1))}
{\Ga_q(\la j +a+1)\Ga_q(\la j +b+1)\Ga_q(\la)}\nonumber\\
&=& \frac{1}{\Ga_{q^{\la}}(n+1)} \frac{1}{(\Ga_q(\la +1))^n}
\prod_{j=0}^{n-1}\frac{\Ga_q(\la j+a+b+1)\Ga_q(\la(j+1)+1)}
{\Ga_q(\la j +a+1)\Ga_q(\la j +b+1)} \label{red-mor}
\end{eqnarray}
where
\begin{equation}\label{bert}
F_0(x):= \prod_{i=1}^n (x_i;q)_a\,(\frac{q}{x_i};q)_b\;\prod_{1\leq i
<j\leq n} (\frac{x_j}{x_i};q)_{\la}\,(q\frac{x_i}{x_j};q)_{\la-1}
\end{equation}
A lemma of Stembridge \cite{stem} gives that (3.4) is equivalent to the
original $q$-Morris identity i.e.~(2.15) with $N_1 = 0$, $N_0 = n$.
Zeilberger's proof of the ``reduced'' $q$-Morris identity relies
on the function $F_0(x)$ being almost anti-symmetric. Indeed,
$F_0(x)=x^{-\de}G_0(x)$ where $\de:=(n-1,n-2,\ldots,2,1,0)$, and
$G_0(x)$ {\it is} anti-symmetric. Thus, the constant term of the function
$F_0(x)$ is just $[x^{\de}]G_0(x)$. As part of his inductive proof, he
essentially uses the equation
\begin{equation}\label{tone}
q^{\beta_n}\left[x^{\beta+\de}(1-qux_n)\prod_{i=1}^{n-1} (1-t\frac{x_n}
{x_i}) \right] G_0(x) = t^{n-1}\left[ x^{\beta+\de}(s-x_n)
\prod_{i=1}^{n-1} (1-t^{-1}\frac{x_n}{x_i}) \right]G_0(x),
\end{equation}
where $u:=q^b$, $s:=q^a$, $t:=q^{\lambda}$,
to relate $[x^{\beta+\de}]G_0(x)$ to $[x^{\de}]G_0(x)$ for various
special values of $\beta$.
This is done by using the anti-symmetry of $G_0(x)$ and his ``Crucial
lemma''

\vspace{2mm}\noindent
{\bf Lemma 3.3}\cite{zeil}\quad If $G_0(x_1,\ldots,x_n)$ is an anti-symmetric
Laurent polynomial, $\gamma\in\IZ^n$, and $\sigma$ a permutation
then $[x^{\sigma(\gamma)}]G_0(x) = {\rm sgn}\:\sigma\;[x^{\gamma}]G_0(x)$.
In particular, if any two components of $\gamma$ are equal, then
$[x^{\gamma}]G_0(x) = 0$.

\vspace{2mm}
As an example of how this is done, let us give a result we shall use
subsequently.

\vspace{2mm}\noindent
{\bf Lemma 3.4}\quad Let $\al_1=(1,0,\ldots,0)$ and $\al_2=(1,0,\ldots,0,-1)$.
Then
\begin{equation}\label{tired.1}
[x^{\al_1+\de}]G_0 = \frac{(s-1)(1-t^n)}{(1-qut^{n-1})(1-t)}\:
[x^{\de}]G_0
\end{equation}
\begin{eqnarray}\label{tired.2}
[x^{\al_2+\de}]G_0 &=& \frac{(t-qs)(1-t^{n-1})}{(1-qst^{n-1})(1-t)}\:
[x^{\de}]G_0  
+ \frac{q(u-t)(1-t^{n-1})}{(1-qst^{n-1})(1-t)}\: [x^{\al_1+\de}]G_0
\end{eqnarray}

\vspace{2mm}\noindent
{\bf Proof}\quad Bearing in mind (\ref{tone}), first look at the expansion
\begin{equation}\label{expand.1}
\prod_{i=1}^{n-1} (1-z\frac{x_n}{x_i}) = \sum_T (-zx_n)^{|T|} x^{-T}
\end{equation}
where the sum is over all $T\subseteq\{1,2,\ldots,n-1\}$, and
$x^{-T}:=\prod_{i\in T}x_i^{-1}$. For each $m$ with $0\leq m\leq n-1$,
there exists a unique set $T$, such that $|T|=m$, and
$x^{\de}x_n^{|T|} x^{-T}$ has distinct exponents; namely $T=\{n-m,
n-m+1,\ldots,n-1\}$. In fact $x^{\de}x_n^{|T|} x^{-T}=x^{\sigma(\de)}$
where $\sigma$ is a permutation with ${\rm sgn}\:\sigma = (-1)^m$. Thus
\begin{equation}\label{today.1}
\left[x^{\de}\prod_{i=1}^{n-1}\left(1-z\frac{x_n}{x_i}\right)\right]G_0
= \sum_{m=0}^{n-1}(-z)^m(-1)^m\,[x^{\de}]G_0 = \left(\frac{1-z^n}{1-z}
\right) [x^{\de}]G_0
\end{equation}
We must also expand
\begin{equation}\label{expand.2}
x_n\prod_{i=1}^{n-1} (1-z\frac{x_n}{x_i}) = \sum_T (-z)^{|T|}(x_n)^{|T|+1} 
x^{-T}
\end{equation}
In this case, there is only one set $T$ such that $x^{\de}x_n^{|T|}x^{-T}$ has
distinct exponents: $T=\{1,2,\ldots,n-1\}$. Moreover, for this set $T$,
$$
x^{\de}x_n^{|T|} x^{-T}=x_n^nx_1^{n-2}x_2^{n-3}\cdots x_{n-2}^1x_{n-1}^0
=x^{\si(\al_1+\de)}
$$
where ${\rm sgn}\:\si=(-1)^{n-1}$. Thus
\begin{equation}\label{today.2}
\left[x^{\de}x_n\prod_{i=1}^{n-1}\left(1-z\frac{x_n}{x_i}\right)\right]G_0
=z^{n-1}[x^{\al_1+\de}]G_0
\end{equation}
If one now uses (\ref{today.1}) and (\ref{today.2}) with $z=t$, $t^{-1}$
in (\ref{tone}), and sets $\beta=0$, the stated result
(\ref{tired.1}) follows.

To prove (\ref{tired.2}), note that the sets $T$ such that
$x^{\al_2+\de}x_n^{|T|} x^{-T}$ has distinct exponents are of the form
$$
T = \left\{ \begin{array}{cc} \emptyset & \\
\{1,n-m+1,n-m+2,\ldots,n-1\} & 1\leq m \leq n-1 
\end{array} \right.
$$
in which case
$$
x^{\al_2+\de}x_n^n x^{-T} =  \left\{ \begin{array}{cc} 
x^{\al_2+\de} & \\
x^{\si(\de)} & {\rm sgn}\:\si = (-1)^{m-1}
\end{array} \right.
$$
Thus
\begin{equation}\label{today.3}
\left[x^{\al_2+\de}\prod_{i=1}^{n-1}\left(1-z\frac{x_n}{x_i}\right)\right]G_0
= [x^{\al_2+\de}]G_0 - z\left(\frac{1-z^{n-1}}{1-z}\right) [x^{\de}]G_0
\end{equation}
Similarly, the sets $T$ such that $x^{\al_2+\de}x_n^{|T|+1}x^{-T}$ has
distinct exponents are of the form $T=\{n-m,n-m+1,\ldots,n-1\}$, 
$0\leq m \leq n-2$, in which case $x^{\al_2+\de}x_n^{|T|+1}x^{-T}=
x^{\si(\al_1+\de)}$, with ${\rm sgn}\:\si=(-1)^m$. Hence
\begin{equation}\label{today.4}
\left[x^{\al_2+\de}x_n\prod_{i=1}^{n-1}\left(1-z\frac{x_n}{x_i}\right)
\right]G_0 = \left(\frac{1-z^{n-1}}{1-z}\right) [x^{\al_1+\de}]G_0
\end{equation}
Again, using (\ref{today.3}), (\ref{today.4}) in (\ref{tone}) (setting
$\beta=\al_2$), with $z=t$, $t^{-1}$ yields (\ref{tired.2}).

\vspace{2mm}
Returning to the proof of the $N_1=2$ case of Conjecture 2.1, 
we first make the substitutions $w_i\rightarrow w_{N_0+1-i}$,
$z_i\rightarrow z_{N_1+1-i}$ (which has no effect on the constant term)
and then follow the arguments in ref \cite{stem}, whereby we replace
$(qz_i/z_j;q)_{\la+1}\rightarrow(qz_i/z_j;q)_{\la}$ and
$(qw_i/w_j;q)_{\la}\rightarrow(qw_i/w_j;q)_{\la-1}$, to obtain an
alternative statement of Conjecture 2.1 in the case $p=2$ 
which reads as follows:
\begin{eqnarray}
{\rm CT}\, \prod_{i=1}^{N_0} (w_i;q)_a (q/w_i;q)_b \prod_{j=1}^{N_1}
(z_j;q)_a (q/z_j;q)_b\;\prod_{1\leq i<j\leq N_0}(\frac{w_j}{w_i};q)_{\la}
(q\frac{w_i}{w_j};q)_{\la-1} \nonumber\\
\times\prod_{1\leq i<j\leq N_1}(\frac{z_j}{z_i};q)_{\la+1}
(q\frac{z_i}{z_j};q)_{\la} \prod_{i=1}^{N_0}\prod_{j=1}^{N_1}
(q\frac{w_i}{z_j};q)_{\la}(\frac{z_j}{w_i};q)_{\la} \nonumber\\
=\frac{1}{\Ga_{q^{\la}}(N_0+1) \Ga_{q^{\la+1}}(N_1+1)}\, D_1(N_1;N_0;
a,b,\la) \hspace{4cm}
\end{eqnarray}
where $D_1(N_1;N_0;a,b,\la)$ is given in (2.15).
In the particular case of $N_1=2$, the function appearing on the
left-hand side of the above equation, call it $F_1(\{z_j\};\{w_j\})$ 
say, is simply related to the $n=N_0+2$ variable function $F_0(x)$ 
in (\ref{bert}). Thus, letting $x_i=w_i$, $1\leq i\leq N_0$, and
$x_{N_0+1}=z_1$, $x_{N_0+2}=z_2$, we have
$$
F_1(x):= (1-t\frac{x_{N_0+1}}{x_{N_0+2}})(1-t\frac{x_{N_0+2}}{x_{N_0+1}})
\prod_{i=1}^{N_0}
(1-t\frac{x_i}{x_{N_0+1}})(1-t\frac{x_i}{x_{N_0+2}})\,F_0(x)
$$
Using the ``reduced'' $q$-Morris identity (\ref{red-mor}), it 
suffices to prove
\begin{equation}\label{g.1}
[x^0]F_1(x) = \frac{(1-t^{N_0+1})(1-t^{N_0+2})(1-qust^{N_0+1})
(1-qt^{N_0+2})}{(1-t)^2(1-qut^{N_0+1})(1-qst^{N_0+1})}\;
[x^0]F_0(x)
\end{equation}
Note that we can rewrite $[x^0]F_1(x)$ in the following form,
\begin{equation}\label{guff.1}
[x^0]F_1(x) = \left[\left((1+t^2)-t\frac{x_{N_0+1}}{x_{N_0+2}}
-t\frac{x_{N_0+2}}{x_{N_0+1}}\right)\prod_{i=1}^{N_0}\left(
1-t\frac{x_{N_0+1}}{x_i}\right)\left(1-t\frac{x_{N_0+2}}{x_i}\right)
x^{\de}\right]\,G_0(x)
\end{equation}
where $\de=(N_0+1,N_0,\ldots,1,0)$, and $G_0(x)$ is anti-symmetric.
Let us now show that each of the terms
$[A(x)\prod_{i=1}^{N_0}(1-tx_{N_0+1}/x_i)(1-tx_{N_0+2}/x_i)x^{\de}]G_0$,
for $A(x)=1+t^2$, $-tx_{N_0+1}/x_{N_0+2}$ and $-tx_{N_0+2}/x_{N_0+1}$
can be expressed in terms of $[x^{\de}]G_0$ and $[x^{\al_2+\de}]G_0$
using the above techniques.

\vspace{2mm}\noindent
{\bf Lemma 3.5}\quad We have 
\begin{eqnarray*}
\left[\prod_{i=1}^{N_0}\left(1-t\frac{x_{N_0+1}}{x_i}\right)
\left(1-t\frac{x_{N_0+2}}{x_i}\right) x^{\de}\right]G_0
&=& B_{N_0}(t) \\ 
\left[-t\frac{x_{N_0+2}}{x_{N_0+1}}\prod_{i=1}^{N_0}\left(1-t
\frac{x_{N_0+1}}{x_i}\right) \left(1-t\frac{x_{N_0+2}}{x_i}\right) 
x^{\de}\right]G_0
&=& tB_{N_0}(t) 
\end{eqnarray*}
where
$$
B_{N_0}(t) = \frac{1}{(1-t)}\left( \frac{(1-t^{N_0+1})}{(1-t)} - 
t\frac{(1-t^{2N_0+2})}{(1-t^2)} \right)[x^{\de}]G_0
$$

\vspace{2mm}\noindent
{\bf Proof}\quad We prove only the first formula, as the proof of
the second is similar. First, expand 
$$
\prod_{i=1}^{N_0}\left(1-t\frac{x_{N_0+1}}{x_i}\right)
\left(1-t\frac{x_{N_0+2}}{x_i}\right)  
=\sum_{n,m,\ga}(-t)^{n+m}c_{n,m,\ga}\,f_{n,m,\ga}(x^{-1})\,
x_{N_0+1}^{n} x_{N_0+2}^{m}
$$
where $f_{n,m,\ga}$ is the monomial $x_1^{\ga_1}\cdots x_{N_0}^{\ga_{N_0}}$ 
with exponents $\ga_i=0$, $-1$ or $-2$, and $c_{n,m,\ga}$ is a positive 
integer. The only terms in this expansion which have distinct exponents when
multiplied by $x^{\de}$ occur when $n\geq m$. Moreover
$$
f_{n,m,\ga}=x_{N_0+1-n}^{-1}x_{N_0+2-n}^{-1}\cdots x_{N_0-m}^{-1}
x_{N_0+1-m}^{-2}\cdots x_{N_0}^{-2},
$$
$c_{n,m,\ga}=1$ and 
$$
f_{n,m,\ga}(x^{-1})x_{N_0+1}^nx_{N_0+2}^mx^{\de}\!=\!
x_1^{N_0+1}\!\cdots x_{N_0-n}^{n+2}
x_{N_0+1}^{n+1}x_{N_0+1-n}^n\!\cdots\!x_{N_0-m}^{m+1}x_{N_0+2}^m
x_{N_0+1-m}^{m-1}\!\cdots\! x_{N_0}^0 \!=\! x^{\si(\de)}
$$
where ${\rm sgn}\,\si=(-1)^{n+m}$. Thus, 
$$
\left[\prod_{i=1}^{N_0}\left(1-t\frac{x_{N_0+1}}{x_i}\right)
\left(1-t\frac{x_{N_0+2}}{x_i}\right) x^{\de}\right]\,G_0
=\sum_{n=0}^{N_0}\sum_{m=0}^n (-1)^{n+m}(-t)^{n+m} [x^{\de}]G_0
$$
which yields the result upon carrying out the summation. 

\vspace{2mm}\noindent
{\bf Lemma 3.6}\quad We have
\begin{eqnarray}
\left[-t\frac{x_{N_0+1}}{x_{N_0+2}}\prod_{i=1}^{N_0}\left(1-t
\frac{x_{N_0+1}}{x_i}\right)
\left(1-t\frac{x_{N_0+2}}{x_i}\right) x^{\de}\right]\,G_0(x)
= -t^{N_0+1}[x^{\al_2+\de}]\,G_0(x) \nonumber\\
- \frac{t^2}{1-t}\left(\frac{1-t^{2N_0}}{1-t^2} - t^{N_0}
\frac{1-t^{N_0}}{1-t} \right) [x^{\de}]\,G_0(x)\label{blah}
\end{eqnarray}

\vspace{2mm}\noindent
{\bf Proof}\quad Again, expand
$$
-t\frac{x_{N_0+1}}{x_{N_0+2}}\prod_{i=1}^{N_0}\left(1-t
\frac{x_{N_0+1}}{x_i}\right) \left(1-t\frac{x_{N_0+2}}{x_i}\right)
=\sum_{n,m,\ga}(-t)^{n+m+1}c_{n,m,\ga}\,f_{n,m,\ga}(x^{-1})\,
x_{N_0+1}^{n+1} x_{N_0+2}^{m-1}
$$
Once more, $f_{n,m,\ga}$ is a monomial in $x_1^{-1},\ldots,x_{N_0}^{-1}$
with exponents no greater than $-2$ and $c_{n,m,\ga}$ is a positive integer. 
The only terms in this expansion which, when
multiplied by $x^{\de}$, have distinct exponents occur when {\it either}
$1\leq m\leq N_0$ and $m-1\leq n\leq N_0-1$ {\it or} $n=N_0$, $m=0$.

In the latter case, $c_{N_0,0,\ga}=1$, 
$f_{N_0,0,\ga}=(x_1x_2\cdots x_{N_0})^{-1}$ and
\begin{equation}\label{poss.0}
f_{N_0,0,\ga}\,x_{N_0+1}^{N_0+1}x_{N_0+2}^{-1}x^{\de} = x^{\si(\al_2+\de)}
\hspace{2cm} {\rm sgn}\,\si=(-1)^{N_0}
\end{equation}

In the former case the monomials $f_{n,m,\ga}$ take one of the 
$n-m+1$ possible forms
\begin{equation}\label{poss.1}
f_{n,m,\ga} = \left\{\begin{array}{c} 
x_{N_0-n}^{-2}x_{N_0-n+2}^{-1}x_{N_0-n+3}^{-1}\cdots
x_{N_0-m+1}^{-1}x_{N_0-m+2}^{-2}\cdots x_{N_0}^{-2} \\[2mm]
x_{N_0-n}^{-1}x_{N_0-n+1}^{-2}x_{N_0-n+3}^{-1}\cdots
x_{N_0-m+1}^{-1}x_{N_0-m+2}^{-2}\cdots x_{N_0}^{-2} \\[2mm]
\vdots \\
x_{N_0-n}^{-1}x_{N_0-n+1}^{-1}\cdots x_{N_0-m+1}^{-1}
x_{N_0-m}^{-2}x_{N_0-m+2}^{-2}\cdots x_{N_0}^{-2} 
\end{array}\right.
\end{equation}
as well as the additional form
\begin{equation}\label{poss.2}
f_{n,m,\ga} = x_{N_0-n}^{-1}x_{N_0-n+1}^{-1}\cdots x_{N_0-m}^{-1}
x_{N_0-m+1}^{-1}x_{N_0-m+2}^{-2}\cdots x_{N_0}^{-2} 
\end{equation}
For the monomials (\ref{poss.1}) the corresponding $c_{n,m,\ga}=1$,
but for the monomial (\ref{poss.2}), $c_{n,m,\ga}=n-m+2$. Moreover,
in the former case for each $f_{n,m,\ga}$, we have
$f_{n,m,\ga}(x^{-1})x_{N_0+1}^{n+1}x_{N_0+2}^{m-1}
x^{\de}=x^{\si(\de)}$, where ${\rm sgn}\,\si=(-1)^{n+m+1}$, while for
the latter, ${\rm sgn}\,\si=(-1)^{n+m}$. Hence, combining contributions
from (\ref{poss.1}), (\ref{poss.2}) and (\ref{poss.0}), we get
\begin{eqnarray*}
\mbox{LHS of (\ref{blah})} &=& \sum_{m=1}^{N_0}\sum_{n=m-1}^{N_0-1}
(-t)^{n+m+1}\left((n-m+1)(-1)^{n\!+\!m\!+\!1} +(n\!-\!m+\!2)(-1)^{n+m}\right)
[x^{\de}]G_0 \\
&&-t^{N_0+1}[x^{\al_2+\de}]G_0 
\end{eqnarray*}
which produces the required result after summation.

\vspace{2mm}
Lemmas 3.5 and 3.6 show that everything on the right hand side of
(\ref{guff.1}) can be expressed in terms of $[x^{\de}]G_0$ and 
$[x^{\al_2+\de}]G_0$. However, by eliminating $[x^{\al_1+\de}]G_0$
from (\ref{tired.1}) and (\ref{tired.2}), we have
$$
[x^{\al_2+\de}]G_0 = \frac{1}{(1-qst^{N_0+1})}\left(\frac{1-t^{N_0+1}}
{1-t}\right)\left\{(t-qs)+q(u-t)\frac{(s-1)(1-t^{N_0+2})}{(1-qut^{N_0+1})
(1-t)} \right\}\;[x^{\de}]G_0
$$
Thus, from (\ref{guff.1})
\begin{eqnarray*}
\frac{[x^0]F_1}{[x^{\de}]G_0} = (1+t+t^2)\times\frac{1}{1-t}
\left( \frac{(1-t^{N_0+1})}{(1-t)} - t\frac{(1-t^{2N_0+2})}
{(1-t^2)} \right) \\
- \frac{t}{1-t} \left( t\frac{(1-t^{2N_0})}{(1-t^2)} - t^{N_0+1}\frac{
(1-t^{N_0})}{(1-t)} \right) \\
-t\frac{1}{(1-qst^{N_0+1})}\left(\frac{1-t^{N_0+1}}
{1-t}\right)\left\{(t-qs)+q(u-t)\frac{(s-1)(1-t^{N_0+2})}{(1-qut^{N_0+1})
(1-t)} \right\}
\end{eqnarray*}
Simplification of this expression yields the desired result (\ref{g.1}).

\vspace{.5cm}
\noindent
{\bf 4. CONJECTURES INVOLVING MACDONALD POLYNOMIALS}

\vspace{.3cm}
\noindent
{\bf 4.1 A Gram-Schmidt construction}

\vspace{.2cm}
\noindent
The $q$-generalization (2.10) of (1.1) can be used as a weight function in the
Gram-Schmidt construction of the $q$-generalization of the polynomials
satisfying conditons (i) and (ii) of Section 1. Thus define the
$q$-generalization of the inner product (1.2) by
\renewcommand{\theequation}{4.1}
\begin{eqnarray}
\lefteqn{
\langle f|g \rangle_{N_0,\dots,N_p;\lambda;q}} \nonumber \\& & := 
  \prod_{l=1}^{N_0} \int_{-1/2}^{1/2} d y_l \,
\prod_{\alpha = 1}^p \prod_{l=1}^{N_\alpha}
 \int_{-1/2}^{1/2} d x_l^{(\alpha)} \,|\psi_0(\{z_j^{(\alpha)}
\}_{\alpha = 1,\dots,p \atop j=1,\dots,N_\alpha},
\{w_j\}_{j=1,\dots,N_0};q)|^2 f^* \, g.
\end{eqnarray}
and define the $q$-generalization of the symmetric polynomials $p_\kappa$,
 $p_\kappa(w_1, \dots,w_{N_0};q)$ say, by properties (i) and (ii) of Section 1
with the inner product therein replaced by (4.1).
Note from condition (i) that
$$
p_{1^k}(w_1, \dots,w_{N_0};q) =   m_{1^k} := \sum_{1 \le j_1 < \dots < j_k \le
N_0} w_{j_1} w_{j_2} \dots w_{j_k}.
\eqno (4.2)
$$

Based on some exact computer generated data, and the conjecture (1.3) for the
$q=1$ case, we make the following conjecture.

\vspace{.2cm}
\noindent
{\bf Conjecture 4.1} \quad We have
$$
p_\kappa(w_1, \dots,w_{N_0};q) = P_\kappa(w_1, \dots,w_{N_0};qt^p,t), \qquad {\rm
where}
\quad t := q^\lambda
\eqno (4.3)
$$
and $P_\kappa$ denotes the Macdonald polynomial. In the limit $q \rightarrow
1$, Conjecture 4.1 reduces to (1.5). Note in particular that this 
agrees with the original conjecture \cite[Conj.~2.4]{forr2} only in the 
case $p=1$.
\vspace{.2cm}

We have obtained exact computer generated data for two further conjectures
which are closely related to this result. In relation to the first conjecture
we note from Conjecture 4.1 that $p_\kappa$ is independent of $N_1, \dots, N_p$.
Analogous to the $q=1$ case \cite[Conj.~2.3]{forr2}, this can be understood in
terms of a  conjecture which generalizes Conjecture 2.2.

\vspace{.2cm}
\noindent
{\bf Conjecture 4.2} \quad Let $h = h(w_1, \dots, w_{N_0})$ be a Laurent
polynomial of the form $h = \sum_{\bar{\sigma} \le \bar{\rho}} c_\sigma
m_\sigma$, where $\rho = (\rho_1, \dots, \rho_{N_0})$, $|\rho_1| \ge \dots \ge 
|\rho_{N_{0}}|$ and $\bar{\rho} = (|\rho_1|, \dots, |\rho_{N_0}|)$, and let
$$
D_p(N_1, \dots,N_p;N_0;a,b,\lambda;q)[h] := 
\langle h | A(\{ z_j^{(\alpha)}\}_{\alpha = 1,\dots,p \atop j=1,\dots,N_\alpha},
 \{ w_j\}_{j=1,\dots,N_0};q) \rangle_{N_0,\dots,N_p;\lambda;q}
$$
where $ A(\{ z_j^{(\alpha)}\},
 \{ w_j\;q)$ is given by (2.16b). For $N_p \ge N_j - 1$
 and $N_j \ge |\rho_1| - 1$  $(j=1,\dots,p-1)$ we conjecture that
$$
{D_p(N_1, \dots,N_{p-1},N_p + 1;N_0;a,b,\lambda;q)[h] \over
D_p(N_1, \dots,N_{p-1},N_p;N_0;a,b,\lambda;q)[h]}
$$
is given by the r.h.s.~of Conjecture 2.2.

\vspace{.2cm}
{}From Conjecture 4.2 it follows that for $N_p \ge N_j - 1 $ and
 $N_j \ge |\rho_1| - 1$ $(j=1,\dots,p-1)$
$$
D_p(N_1, \dots,N_p;N_0;0,0,\lambda;q)[h] = f_{p-1}(N_1,\dots,N_{p-1};N_0;
\lambda;q)[h]A_{p,|\rho_1|}(N_1,\dots,N_p;N_0;\lambda;q)
\eqno (4.4)
$$
where $A_{p,|\rho_1|}$ denotes the r.h.s.~of (2.17) with $N_p$ replaced by
$l$ and the product formed over $l$ from $l=|\rho_1|-1$ to $N_p-1$. From the
symmetry
$$
D_p(N_1,\dots,N_{p-2},N,N-1;N_0;0,0,\lambda;q)[h] =
D_p(N_1,\dots,N_{p-2},N-1,N;N_0;0,0,\lambda;q)[h]
$$
we see from (4.4) that for some function $f_{p-2}$ which is independent of
$N_p$ and $N_{p-1}$
$$
f_{p-1}(N_1,\dots,N_{p-1};N_0;
\lambda;q)[h] = A_{p-1,|\rho_1|}(N_1,\dots,N_{p-1};N_0;\lambda;q)
f_{p-2}(N_1,\dots,N_{p-2};N_0;
\lambda;q)[h]
$$
where
$$
A_{p-1,|\rho_1|}(N_1,\dots,N_{p-1};N_0;\lambda;q) :=
\prod_{n=|\rho_1|}^{N_{p-1}}
{A_{p,|\rho_1|}(N_1,\dots,N_{p-2},n-1,n;N_0;\lambda;q) \over
A_{p,|\rho_1|}(N_1,\dots,N_{p-2},n,n-1;N_0;\lambda;q)}.
$$
Thus
\begin{eqnarray*}\lefteqn{
D_p(N_1, \dots,N_p;N_0;0,0,\lambda;q)[h] = } \\
& &  f_{p-2}(N_1,\dots,N_{p-2};N_0;
\lambda;q)[h]A_{p,|\rho_1|}(N_1,\dots,N_p;N_0;\lambda;q)
A_{p-1,|\rho_1|}(N_1,\dots,N_{p-1};N_0;\lambda;q).
\end{eqnarray*}
Proceeding in this fashion we see that the dependence of $D_p$ on $N_1,
\dots,N_p$ factorizes from the dependence on $h$, $N_0$ and $p$
and thus cancels out of the ratio of inner products which define the
coefficients in the Gram-Schmidt procedure.

In relation to the second conjecture, we recall \cite{macd} that the Macdonald
polynomial \\ $P_\kappa(w_1, \dots,w_{N_0};q,t)$ is an eigenfunction of the
(mutually commuting) operators
$$
M_{N_0}^{(r)}(q,t) := \sum_{I} A_I(w;t)\,\prod_{i\in I} T_{q,w_i}
$$
summed over all $r$-element subsets of $\{1,2,\ldots,N_0\}$,
where 
$$
A_I(w;t) = t^{r(r-1)/2} \left ( \prod_{i\in I \atop j \notin I}
{t w_i - w_j \over w_i - w_j} \right ) 
\eqno (4.5a)
$$
and $ T_{q,w_i}$ is the $q$-shift operator with action
$$
 T_{q,w_i} f(w_1, \dots,w_{i-1},w_i,w_{i+1}, \dots,w_{N_0}) =
 f(w_1, \dots,w_{i-1},qw_i,w_{i+1}, \dots,w_{N_0}).
\eqno (4.5b)
$$
Furthermore, the fact that the set of Macdonald polynomials $\{P_\kappa(w_1,
\dots, w_{N_0};q,q^\lambda\}_{\kappa}$ are orthogonal with respect to the inner
product (4.1) with $p=0$ follows from the eigenfunction property and the fact
that $M_{N_0}^{(r)}(q,t)$ is Hermitian with respect to this inner product:
$$
\langle M_{N_0}^{(r)}(q,t)f|g\rangle_{N_0;\lambda;q} = 
\langle f| M_{N_0}^{(r)}(q,t) g\rangle_{N_0;\lambda;q}
\eqno (4.6)
$$

For general $p$ we have obtained exact computer-generated data which suggest
a result similar to (4.6).

\vspace{.2cm}
\noindent
{\bf Conjecture 4.3} \quad Let $m_\kappa$ and $m_\mu$ be monomial symmetric
functions and suppose 
$$
N_j \ge \left \{ \begin{array}{cc}{\rm min}(\kappa_1,\mu_1), & \kappa_1
\ne \mu_1 \\ \kappa_1 - 1, & \kappa_1 = \mu_1 \end{array} \right .
$$ 
$(j=1,\dots,p)$.
We have
$$
\langle M_{N_0}^{(r)}(qt^p,t) m_\kappa|m_\mu\rangle_{N_0,\dots,N_p;\lambda;q}
= \langle m_\kappa| M_{N_0}^{(r)}(qt^p,t) m_\mu\rangle_{N_0,\dots,N_p;\lambda;q}.
$$

\vspace{.2cm}\noindent
{\it Remarks}\\
1.\quad Since the set of functions $\{P_\kappa(w_1,
\dots, w_{N_0};qt^p,t\}_{\kappa}$ are eigenfunctions of $ M_{N_0}^{(r)}
(qt^p,t)$ with distinct eigenvalues, we see that Conjecture 4.1 follows 
as a corollary of Conjecture 4.3.\\[2mm]
2.\quad
As $q\rightarrow 1$ Conjecture 4.3 implies  that
the Laplace-Beltrami operator
$$
D_{N_0}(p+1/\la) = -\frac{(p+1/\la)^2}{2} \sum_{i=1}^{N_0} w_i\,
\frac{\partial^2}{\partial w_i^2} - (p+1/\la) \sum_{i\neq j}
\frac{w_i^2}{w_i-w_j}\,\frac{\partial}{\partial w_i},
$$
which is (up to a shift by a constant) the limit $q\rightarrow 1$ of
$$
\frac{1}{(q^{\la}-1)^2}\left(M_{N_0}^{(2)} -(N_0-1)\,M_{N_0}^{(1)}
+\frac{1}{2}N_0(N_0-1) \,M_{N_0}^{(0)} \right), 
$$
is Hermitian with repect to the multi-component Jack inner product
(1.2), given the same constraints on $N_1,\ldots,N_p$.

\vspace{.3cm}
\noindent
{\bf 4.2 Normalization integral}

\vspace{.2cm}
\noindent
In the $q = 1$ case the normalization of (1.3) with respect to the inner product
(1.2) has been conjectured in ref.[2, Section 3]. For general $q$ the
corresponding normalization is given by the inner product
$$
\langle  P_\kappa(w_1, \dots,w_{N_0};qt,t) |
 P_\kappa(w_1, \dots,w_{N_0};qt,t) \rangle_{N_0,\dots,N_p;\lambda;q}
=: {\cal N}_p^\kappa (N_1, \dots,N_p;N_0;\lambda;q)
\eqno (4.7)
$$
with $\kappa_1$ restricted as in condition (i) of Section 1.
As a corollary to Conjecture 4.2, the dependence on $N_1, \dots, N_p$
factorizes from the dependence on $\kappa$, $N_0$ and $p$, so it suffices to
consider the case $N_1=\dots=N_p=\kappa_1$. Explictly, in the case $p=1$
Conjecture 4.2 gives
$$
{\cal N}_1^\kappa(N_1;N_0;\lambda;q)=
{\cal N}_1^\kappa(\kappa_1;N_0;\lambda;q)
\prod_{j = \kappa_1}^{N_1-1} {[(\lambda + 1)(j+1)]_q \over \Gamma_q (1 +\lambda)}
\Big ( (\lambda + 1)j + \lambda N_0 + 1;q \Big )_\lambda
\eqno (4.8)
$$

To obtain a conjecture for the evaluation of ${\cal
N}_1^\kappa(N_1;N_0;\lambda;q)$ we have obtained some exact computer generated
data, as well as an analytic result in the case $\kappa = 1^k$, $\lambda = 1$
(see the next section). These results, and the corresponding conjecture 
\cite[Conjecture 3.1]{forr2} in the $q=1$ case, suggest a closed form expression
for general
$\kappa$.

\vspace{.2cm}
\noindent
{\bf Conjecture 4.4} \quad Let $f_j$ denote the frequency of the integer $j$ in
the partition $\kappa$ so that $\kappa = \kappa_1^{f_{\kappa_1}}
 (\kappa_1 - 1)^{f_{\kappa_1 - 1}} \dots 1^{f_1}$. We have
\begin{eqnarray*}
\lefteqn{{\cal N}_1^{f_{\kappa_1} f_{\kappa_1 -1}\dots f_1}
(\kappa_1;N_0;\lambda;q)
={[\kappa_1]_q! \Gamma_{q^\lambda}(\lambda N_0 + 1) \over 
(\Gamma_q(1+\lambda))^{N_0 + \kappa_1}}
 \prod_{j=1}^{\kappa_1}\Big (\lambda f_j + 1;q \Big )_\lambda  }\\
& & 
 \times
 \prod_{j=1}^{\kappa_1}{\Big ((\lambda + 1)j + 1 + \lambda
{\displaystyle \sum_{k=1}^{\kappa_1 + 1 -j}} f_{j+k-1};q \Big )_\lambda
\Big ((\lambda + 1)(j-1) + 1 + \lambda
(N_0 - {\displaystyle \sum_{k=1}^{\kappa_1 + 1 -j}} f_{j+k-1});q 
\Big )_\lambda \over
\Big ( (\lambda + 1)j + 1 + \lambda f_j;q \Big )_\lambda} \: \: .
\end{eqnarray*}

\vspace{.5cm}
\noindent
{\bf 5. AN EXPLICIT GRAM-SCHMIDT CONSTRUCTION AND COMPUTATION OF THE
NORMALIZATION}

\vspace{.3cm}
\noindent
{\bf 5.1 A $q$-determinant method}

\vspace{.2cm}
\noindent
In ref.~\cite[Proposition 2.1]{forr2} a determinant method was used to prove, for
$q=1$,  Conjecture 4.1 in the case $p = \lambda = 1$, $\kappa = 21^k$ and
Conjecture 4.4 in the case $p = \lambda = 1$. It is possible to $q$-generalize
the determinant method and thus prove these results for general $q$.

\vspace{.2cm}
\noindent
{\bf Proposition 5.1} \quad
For $p = \lambda = 1$ and $N_1 \ge 1$ we have
$$
\langle  m_{1^{k+2}} |  m_{1^{k+2}} \rangle_{N_0;N_1;1;q} =
 [N_0]_q! [N_1]_{q^2}![k+3]_q[N_0 -1 -k]_q \prod_{l =1}^{N_1 - 1}[N_0 + 2N_1 + 1
- 2l]_q
\eqno (5.1)
$$
and
$$
\langle   m_{1^{k+2}}| s_{21^k}\rangle_{N_0;N_1;1;q}  = 
q [N_0]_q! [N_1]_{q^2}![k+1]_q[N_0 -1 -k]_q \prod_{l =1}^{N_1 - 1}[N_0 +
2N_1 + 1 - 2l]_q
\eqno (5.2) 
$$
where
$s_\kappa = s_\kappa (w_1, \dots, w_{N_0})$ denote the Schur polynomial.
Consequently
$$
p_{21^k}(w_1, \dots, w_{N_0};q) = s_{21^k} - q{[k+1]_q \over [k+3]_q}
s_{1^{k+2}}.
\eqno (5.3)
$$

\vspace{.2cm}
\noindent
{\bf Proof} \quad We will first show how to deduce (5.3) from (5.1) and (5.2).
Now, since the Schur polynomials can be written
$$
s_\kappa = m_\kappa + \sum_{\mu < \kappa} a_\mu' m_\mu
$$
for some coefficients $a_\mu'$, the condition (i) (recall Section 1) in the
definition of
$p_\kappa$ can be rewritten as
$$
p_\kappa(w_1, \dots, w_{N_0};q) = s_\kappa +  \sum_{\mu < \kappa} a_\mu s_\mu.
$$
Furthermore,
$$
s_{1^n} = m_{1^n} = p_{1^n}(w_1, \dots, w_{N_0};q),
$$
so by the Gram-Schmidt procedure
$$
p_{21^k}(w_1, \dots, w_{N_0};q) = s_{21^k}(w_1, \dots, w_{N_0}) -
{\langle    m_{1^{k+2}} |  s_{21^{k}} \rangle_{N_0;N_1;1;q} \over
\langle  m_{1^{k+2}} |  m_{1^{k+2}} \rangle_{N_0;N_1;1;q}}
 s_{1^{k+2}}(w_1, \dots,
w_{N_0})
\eqno (5.4)
$$
Substitution of (5.1) and (5.2) into (5.4) gives (5.3).

Next we take up the task of deriving the results (5.1) and (5.2). We first
transform the integrand into a form symmetric in $\{z_j\}$ and $\{w_j\}$ by
appealing to a lemma of Kadell \cite[lemma 4]{kade}, which for any $f$ and
$g$ symmetric
in
$\{z_j\}$ and $\{w_j\}$ gives the identity
$$
\langle  g |  f \rangle_{N_0;N_1;1;q} = 
{[N_0]_q! \over N_0!}{[N_1]_{q^2}! \over N_1!}
\prod_{l=1}^{N_1} \int_{-1/2}^{1/2}dx_l 
\prod_{l=1}^{N_0} \int_{-1/2}^{1/2}dy_l \,
F(\{z_j\},\{w_j\};q) f(\{z_j\},\{w_j\};q)
$$
$$
\hspace*{3cm} \times g(\{z_j^*\},\{w_j^*\};q)
$$
where
$$
F(\{z_j\},\{w_j\};q) =  F_1(\{w_j\}) F_2(\{z_j\},\{w_j\};q)
$$
with
$$
F_1 = \prod_{1 \le j < k \le N_0}(w_k - w_j),
$$
\begin{eqnarray*}\lefteqn{
F_2 = \prod_{1 \le j < k \le N_0}  (w_k^{-1} - w_j^{-1})} \\ && \times
\prod_{j=1}^{N_0} \prod_{\alpha = 1}^{N_1} \Big (1 - {z_\alpha \over w_j}\Big 
)\Big (1 -q {w_j \over z_\alpha}\Big 
)
\prod_{1 \le \alpha < \beta \le N_1}\! \Big (1 -{z_\beta \over
 z_\alpha}\Big )\Big (1 -{z_\alpha \over
 z_\beta}\Big ) \Big (1 -q{z_\beta \over
 z_\alpha}\Big )\Big (1 -q{z_\alpha \over
 z_\beta}\Big )
\end{eqnarray*}
(this factorization of $F$ is chosen for later convenience).
Thus
\renewcommand{\theequation}{5.5}
\begin{eqnarray}\lefteqn{
\langle m_{1^{k+2}}| m_{1^{k+2}} \rangle_{N_0;N_1;1;q} =
{[N_0]_q! \over N_0!}{[N_1]_{q^2} \over N_1!}} \nonumber \\
& & \times \prod_{l=1}^{N_1} \int_{-1/2}^{1/2}dx_l
\prod_{l=1}^{N_0} \int_{-1/2}^{1/2}dy_l F(\{z_j\},\{w_j\};q)
m_{1^{k+2}}(\{w_j\})m_{1^{k+2}}(\{w_j^*\})
\end{eqnarray}
and similarly the inner product in (5.2).

Consider now the task of evaluating the integral in (5.5). Our method 
is to write $F$ in terms of determinants. From the 
Vandermonde determinant identity
$$
 \prod_{1 \le j < k \le N_0 + 2N_1}(u_k - u_j) = \det [u_j^{k-1}]_{j,k = 1,
\dots, N_0 + 2N_1},
$$
with
$$
u_j = w_j^{-1} \: (j=1,\dots,N_0) \quad
u_{j+N_0} = z_j^{-1} \: (j=1,\dots,N_1) \quad
u_{j+N_0+N_1} = qz_j^{-1} \: (j=1,\dots,N_1), 
$$
straightforward manipulation gives
$$
F_2(\{z_j\},\{w_j\};q)  =
(-1)^{N_1 + N_0N_1}(1 - q)^{-N_1} q^{-N_1(N_1 - 1)/2}
\prod_{j=1}^{N_1} z_j^{1 + 2(N_1 - 1) + N_0}
\prod_{j=1}^{N_0} w_j^{N_1} 
$$
$$
\times \det \left [ \begin{array}{l} [w_j^{-(l-1)}]_{j =1, \dots, N_0 \atop
l =1, \dots, N_0 + 2N_1 } \\
\left [ { z_j^{-(l-1)} \atop (qz_j^{-1})^{(l-1)}} \right ]_{j =1, \dots, N_1 \atop
l =1, \dots, N_0 + 2N_1 } \end{array} \right ]
\eqno (5.6)
$$
(the block notation in (5.6) indicates successive rows; thus the row with
elements $z_1^{-(l-1)}$ is followed by $(qz_1^{-1})^{l-1}$, which is
followed by $z_2^{-(l-1)}$ etc.).
Also, since $m_{1^{k+2}} = s_{1^{k+2}}$, from the
 determinant formula for the
Schur polynomials we have
$$
F_1 s_{1^{k+2}} = \det [ w_j^{l + \kappa_{N_0 - j +1} - 1}
]_{j,l = 1, \dots, N_0}
\eqno (5.7)
$$
where $\kappa_j = 1$ $(j=1, \dots, k+2)$, $\kappa_j = 0$ otherwise.

Since (5.6) and (5.7) are antisymmetric with respect to interchanges
of $w_1, \dots, w_{N_0}$, in the integral (5.5) we can replace
(5.7) by $N_0!$ times its diagonal term. In the definition of
$m_{1^{k+2}}(\{w_j^*\})$:
$$
m_{1^{k+2}}(\{w_j^*\}) = \sum_{1 \le j_1 < \dots < j_{k+2} \le
N_0} w_{j_1}^{-1} w_{j_2}^{-1} \dots w_{j_{k+2}}^{-1}
\eqno (5.8)
$$
take the sum outside the integral and multiply all terms from the summand of
(5.8) and the diagonal term of the determinant (5.7) into appropriate rows of
(5.6). 
Row-by-row integration of the determinant with respect to 
$w_1, \dots, w_{N_0}$ gives a non-zero contribution in row $j$
only in column
$$
l = N_1 + j + \kappa_{N_0 - j + 1} - \xi_j, \qquad
\xi_j := \left \{ \begin{array}{ll} 1 \quad &{\rm if} \quad j=j_1, \dots,j_{k+2}
\\ 0  
\quad &{\rm otherwise} \end{array}
\right .
\eqno (5.9a)
$$
and this term is equal to unity.
For these non-zero columns to be distinct and the determinant thus non-zero we
require
$$
\{j_1, \dots, j_{k+2} \} =\{1, \dots, \nu, N_0 - k -1, \dots, N_0    -\nu \}.
\eqno (5.9b)
$$
for some $\nu = 0, \dots, k+2$. Assuming this condition and expanding the
integrated determinant by the non-zero columns gives, after expanding the
remaining terms and grouping in pairs
\renewcommand{\theequation}{5.10a}
\begin{eqnarray}\lefteqn{
\langle  m_{1^{k+2}} | m_{1^{k+2}} \rangle_{N_0;N_1;1;q}}\nonumber \\ & =&
{[ N_0]_q![N_1]_{q^2} \over N_1!}(1-q)^{-N_1} q^{-N_1(N_1 - 1)/2} \sum_{P(2
\alpha) > P(2\alpha - 1)}
\epsilon (P)
\prod_{\alpha = 1}^{N_1} (q^{P(2 \alpha)-1} -q^{ P(2 \alpha - 1)-1})\nonumber \\
& & \times
\int_{-1/2}^{1/2} dx_\alpha \, z_\alpha^{N_0 + 2N_1 + 1 - P(2 \alpha) -
P(2 \alpha - 1)}
\end{eqnarray}
where
$$
P(\alpha) \in \{ 1, \dots, N_1-1\} \cup \{N_1 + \nu \}
\cup \{ N_1 + N_0 - \nu +1 \} \cup \{ N_1 + N_0 +2, \dots, N_0 + 2N_1 \}.
\eqno (5.10b)
$$
A non-zero contribution to (5.9) requires
$$
P(2 \alpha - 1) = N_0 + 2N_1 + 1 - P(2 \alpha)
\eqno (5.11a)
$$
and
$$
P(2 \alpha ) \in \{ N_1 + N_0 +1 - \nu \} \cup \{ N_1 + N_0 +2, \dots, 2N_1 + N_0
\}.
\eqno (5.11b)
$$
Each of the $N_1!$ different choices (5.11b) give the same contribution to
(5.9) and so
\renewcommand{\theequation}{5.12}
\begin{eqnarray}\lefteqn{
\langle  m_{1^{k+2}} |  m_{1^{k+2}} \rangle_{N_0;N_1;1;q}}\nonumber \\& &  = 
[N_0]_q! [N_1]_{q^2}!(1-q)^{-N_1} q^{-N_1(N_1 - 1)/2}
 \prod_{l =1}^{N_1 - 1}(q^{N_0 + 2N_1 - l} - q^{l - 1})
 \sum_{\nu = 0}^{k+2} (q^{N_0 + N_1 - \nu} - q^{ N_1 + \nu - 1})\nonumber \\
\end{eqnarray}
which after straightforward simplification gives (5.1).

The computation of the analogue of (5.5) for the inner product (5.2) is very
similar. In place of (5.7) we have
$$
F_1s_{21^{k}}(\{w_j\}) = \det [w_j^{l+\kappa_{N_0 - j+1} - 1}]_{j,l = 1,\dots,
N_0}
\eqno (5.13)
$$
where $\kappa_1 = 2$, $\kappa_j = 1$ $(j=2,\dots,k+1)$, $\kappa_j = 0$
otherwise. After integration over $w_1, \dots, w_{N_0}$ the condition (5.9a)
still gives the column number of the non-zero entry in row $j$.
 For the non-zero columns to be distinct and the condition (5.11a) to hold
we see that in place of (5.9b) we require
$$
\{ j_1, \dots, j_{k+2} \} = \{1, \dots, \nu, N_0 - k, \dots, N_0 - \nu + 1\}
$$
for some $\nu = 1, \dots, k+1.$ This shows that $\langle  
m_{1^{k+2}} |  s_{21^{k}}
\rangle_{N_0;N_1;1;q}$ is given by the r.h.s.~of (5.10a) with condition (5.10b),
which simplifies down to the r.h.s.~of (5.12), the only difference being that
the summation over $\nu$ is now from $\nu = 1$ to $k+1$. After evaluating
the sum the result (5.2) follows.

\vspace{.2cm}
Proposition 5.1 immediately establishes Conjecture 4.4 in the case
$p=\lambda=1$. To prove Conjecture 4.1 in the case $p=\lambda=1$, $\kappa =
21^k$ it is necessary to identify (5.3) as the corresponding Macdonald
polynomial. In the $q=1$ case this can be done by appealing to a theorem of
Stanley \cite[Prop.~7.2]{stan} which gives the explicit expansion of
$J_{2^p1^q}^{(\alpha)}$ in terms of monomial symmetric functions. As we know of
no corresponding result for the Macdonald polynomials, it remains to show that
$P_{21^k}(w_1, \dots, w_{N_0};q^2,q)$ is given by the r.h.s.~of (5.3).
This can be done by the characterisation of the Macdonald polynomial as an
eigenfunction of the operator (4.5a) with $r=1$.

\vspace{.3cm}
\noindent
{\bf 5.2 Expansion of $P_{21^k}(w_1, \dots, w_{N_0};q^2,q)$ in terms of Schur
polynomials}

\vspace{.2cm}
\noindent
\setcounter{equation}{13}
\renewcommand{\theequation}{5.\arabic{equation}}
We know that $P_{21^{n-2}}\equiv P_{21^{n-2}}(w_1,\ldots,w_{N_0};q,t)$ 
must have the form 
\begin{equation}\label{leave.0}
P_{21^{n-2}} = s_{21^{n-2}} + \ga\,s_{1^n} = m_{21^{n-2}}
+ (\ga + n-1)\,m_{1^n}
\end{equation}
where $\gamma$ is to be determined.
The action of the operator $M^{(1)}_{N_0}(q,t)$ (see (4.4a)), on the
monomial symmetric functions is given explicitly by\cite{macd}
\begin{equation}\label{leave.1}
M^{(1)}_{N_0}(q,t)\,m_{\ka} = \sum_{\al} \sum_{i=1}^{N_0} 
t^{N_0-i}q^{\al_i}\,s_{\al}
\end{equation}
where the outer sum is over all derangements $\al\in\IN^{N_0}$ of
the partition $\ka$. Due to the modification rules for Schur
functions associated with unordered partitions
(if $\ka_i<\ka_{i+1}$ for any $i$, then $s_{(\ldots,\ka_i,\ka_{i+1},\ldots)} 
= - s_{(\ldots,\ka_{i+1}-1,\ka_{i}+1,\ldots)}$; in particular
$s_{(\ldots,\ka_{i},\ka_{i}+1,\ldots)}=0$), the only
distinct permutations $\al$ of $\ka=(2,1^{n-2},0^{N_0-n+1})$ for which 
$s_{\al}$ is non-zero are of the form $(2,1^{n-2},0^{N_0-n+1})$, or
$(1^p,0,2,1^{n-2-p},0^{N_0-n})$, $p=0,1,\ldots,n-2$. It thus follows from
(\ref{leave.1}) that 
\begin{equation}\label{le.3} 
M^{(1)}_{N_0}(q,t)\,m_{21^{n-2}} = A_1\,s_{21^{n-2}} + A_2\,s_{1^n}
\end{equation}
where
\begin{eqnarray*}
A_1 &=& q^2t^{N_0-1} + qt^{N_0-n+1}[n-2]_t + [N_0-n+1]_t \\ 
A_2 &=& -q^2 t^{N_0-n} [n-1]_t -
q\left( (n-2)(t^{N_0-n}+t^{N_0-1}) + (n-3)t^{N_0-n+1}[n-2]_t\right) \\
&& - t^{N_0-n+1}[n-1]_t -(n-1)[N_0-n]_t
\end{eqnarray*}
Recalling that $P_{1^n}(w;q,t)\equiv m_{1^n}(w)$ (and hence $m_{1^n}$
is an eigenfunction of $M^{(1)}_{N_0}(q,t)$) it follows from (\ref{leave.0}) and
(\ref{le.3}) that
$$
M^{(1)}_{N_0}(q,t)\,P_{21^{n-2}} = A_1\,s_{21^{n-2}} + \left( A_2 + (\ga+n-1)
e(1^n)\right)\,s_{1^n}
$$
where $e(\kappa):= \sum_{i=1}^{N_0} t^{N_0-i}q^{\kappa_i}$ is the
eigenvalue of $P_{\kappa}(w;q,t)$ under $M^{(1)}_{N_0}(q,t)$.
However
$$
M^{(1)}_{N_0}(q,t)\,P_{21^{n-2}} = e(21^{n-2})\,P_{21^{n-2}} = e(21^{n-2})\,
(s_{21^{n-2}} +\ga s_{1^n})
$$
Equating the coefficients of $s_{1^n}$ in these two equations
(which is permissible since the set of Schur functions $s_{\ka}$, 
for $\ka$ a proper partition, is linearly independent) yields
$$
\ga \,=\, \frac{A_2 + (n-1)e(1^n)}{e(21^{n-2}) - e(1^n)}
\,=\, \frac{(q-t)}{(1-qt^{n-1})}\,[n-1]_t
$$
In the particular case $(q,t)\rightarrow (q^2,q)$ and $n=k+2$, this reproduces
the r.h.s.~of (5.3) and thus completes the proof of
 Conjecture 4.1 in the case
$p=\lambda=1$, $\kappa = 21^k$.

\vspace{1cm}
\noindent
{\bf Appendix}

\noindent
In this appendix an application of a recently derived extension of the
$q$-Morris constant term identity involving the Macdonald polynomial will be
given. The extension is \cite[Theorem 4]{kane}
\renewcommand{\theequation}{A1}
\begin{eqnarray}
\lefteqn{{\rm CT}\, P_\kappa(t_1, \dots, t_n;q,q^\lambda) \prod_{i = 1}^n
(t_i; q)_a \Big ({q \over t_i};q \Big )_b
\prod_{1 \le i < j \le n} \Big ({t_i \over t_j};q \Big )_\lambda
\Big (q{t_j \over t_i};q \Big )_\lambda} \nonumber \\&= & (-1)^{|\kappa|} q^{
\sum_{i=1}^n\kappa_i (\kappa_i+1)/2} \nonumber \\& & \times 
{(q^\lambda;q^\lambda)_n \over (1 - q^\lambda)^n}
\prod_{1 \le i < j \le n} (q^{\kappa_i - \kappa_j} q^{\lambda (j - i)};
q)_\lambda  \prod_{i=1}^n {(q;q)_{a + b + (n-i)\lambda} \over
(q;q)_{a + (n-i)\lambda + \kappa_i}
(q;q)_{b + (i-1)\lambda - \kappa_i}} \nonumber \\&= &
  q^{(b+1)|\kappa|} D_0(n;a,b,\lambda;q)
P_\kappa(1, q^\lambda, \dots,q^{(n-1)\lambda};q,q^\lambda)
{[-b]_{\kappa;q}^{(1/\lambda)} \over
[a+1+(n-1)\lambda]_{\kappa;q}^{(1/\lambda)}}
\end{eqnarray}
where $a,b,\lambda$ are assumed to be non-negative integers,
$D_0(n;a,b,\lambda;q)$ is given by (2.14) with $N_1 = 0$, $N_0 = n$,
\renewcommand{\theequation}{A2}
\begin{eqnarray}
[x]_{\kappa;q}^{(1/\lambda)} &:=& \prod_{j=1}^n {\Gamma_q(x - \lambda (j-1)
+ \kappa_j) \over \Gamma_q(x - \lambda(j-1))} \\
& = &  \prod_{j=1}^n [x-\lambda (j-1) + \kappa_j - 1]_q \dots
[x - \lambda (j-1)]_q
\end{eqnarray}
and we have used the formula \cite{macd}
$$
P_\kappa(1, q^\lambda, \dots,q^{(n-1)\lambda};q,q^\lambda) =
q^{\lambda  \sum_{i=1}^n(i-1)\kappa_i}
\prod_{1 \le i < j \le n}{ (q^{\kappa_i - \kappa_j} q^{\lambda (j - i)};
q)_\lambda \over  (q^{\lambda (j - i)};q)_\lambda}
$$
and the manipulation
$$
{(q;q)_p \over (q;q)_{p-\kappa_i}} = (-1)^{\kappa_i} q^{p \kappa_i}
q^{-\kappa_i (\kappa_i - 1)/2}[-p+\kappa_i-1]_q \dots [-p]_q(1-q)^{\kappa_i}
$$
with $p = b + (i-1)\lambda$.
Using (A1) we can calculate the expansion of the power sums in terms of
Macdonald polynomials. 
 We will require a simple lemma, which was used in a special case in 
ref.~\cite[Proposition 2]{forr3}.

\vspace{.2cm}
\noindent
{\bf Proposition A1} \quad
Let $f(w_1, \dots, w_n)$ be symmetric in $w_1, \dots, w_n$ $(w_j = e^{2 \pi i
y_j})$, periodic of period 1 in each variable $y_j$ and homogenoeous of
integer order
$k$ $(k \ne 0)$. Let $u_\epsilon(w_l)$ have the small-$\epsilon$ expansion
$$
u_\epsilon(w_l) = 1+ \epsilon a(w_l) + O(\epsilon^2).
$$
We have
$$
\lim_{\epsilon \rightarrow 0}{1 \over \epsilon}  
\prod_{l=1}^n \int_{-1/2}^{1/2} dy_l \, u_\epsilon(w_l)\,f(w_1, \dots, w_n)
= n \int_{-1/2}^{1/2}u_\epsilon(w_1)\, dy_1 \prod_{l=2}^n \int_{-1/2}^{1/2}
dy_l\,f(1,w_2, \dots, w_n).
$$

\vspace{.2cm}
\noindent
{\bf Proof} \quad For small-$\epsilon$
$$
\prod_{l=1}^n a_\epsilon(w_l) \: \sim \: 1 + \epsilon \sum_{l=1}^n a
(w_l),
$$
and thus, since $f$ is assumed homogeneous of non-zero integer order
$$
\prod_{l=1}^n \int_{-1/2}^{1/2} dy_l \, u_\epsilon(w_l)\,f(w_1, \dots, w_n)
\: \sim \: \epsilon \prod_{l=1}^n \int_{-1/2}^{1/2} dy_l \,\sum_{l=1}^n a
(w_l) f(w_1, \dots, w_n).
$$
The stated result now follows by using the assumption that $f$ is symmetric to
replace $\sum_{l=1}^n a(w_l)$ in the integrand by $n a(w_1)$, then using the
assumption that $f$ is periodic to replace $w_j$ by $w_1w_j$ $(j=2, \dots,n)$
and finally the fact that $f$ is homogeneous of order $k$ to write
$$
f(w_1,w_1w_2 \dots,w_1 w_n) = w_1^k f(1,w_2, \dots, w_n).
$$
\vspace{.2cm}

The expansion of the power sums is given by the following result.

\vspace{.2cm}
\noindent
{\bf Proposition A2} \quad For $k \in Z_{> 0}$ we have
$$
\sum_{i=1}^n w_i^k = \sum_{|\kappa| = k} {\alpha_{\kappa;q} \over
\langle P_\kappa | P_\kappa \rangle'}
 P_\kappa(w_1, \dots, w_n;q,q^\lambda)
$$
 where
$$
\alpha_{\kappa;q} = {[|\kappa|]_q \Gamma_q(\kappa_1) \over [n]_{q^\lambda}! }
 {\Gamma_q(\lambda n + 1) \over (\Gamma_q( 
\lambda + 1)^n} P_\kappa
(1,q^\lambda,
\dots q^{(n-1)\lambda}) {{[0]'}^{(1/\lambda)}_{\kappa;q} \over
[1 + (n-1) \lambda]_{\kappa;q}^{(1/\lambda)}}
$$
(the dash on ${[0]'}^{(1/\lambda)}_{\kappa;q}$ means that the $j=1$ term
in its definition (A2) is to be omitted)and
\begin{eqnarray*}
\langle P_\kappa | P_\kappa \rangle'
&: = &{1 \over n!}{\rm CT} \, P_\kappa(t_1, \dots, t_n;q,q^\lambda)
P_\kappa(1/t_1, \dots, 1/t_n;q,q^\lambda)
\prod_{1 \le i < j \le n} \Big ({t_i \over t_j};q \Big )_\lambda
\Big ({t_j \over t_i};q \Big )_\lambda \\
& = &  \prod_{1 \le i < j \le n}
{(q^{\kappa_i - \kappa_j + \lambda(j-i)};q)_\infty
(q^{\kappa_i - \kappa_j +1+ \lambda(j-i)};q)_\infty \over
(q^{\kappa_i - \kappa_j + \lambda(j-i+1)};q)_\infty
(q^{\kappa_i - \kappa_j +1+ \lambda(j-i-1)};q)_\infty}
\end{eqnarray*}
(the final equality is given in ref.~\cite{macd}; see also \cite{cher}).

\vspace{.2cm}
\noindent
{\bf Proof}

\vspace{.2cm}
\noindent
First we write the l.h.s.~of (A1) in symmetric form using the lemma of Kadell
used in the proof of Proposition 5.1, and then extend its validity to general
$a,b,\lambda$ by using the integral (2.9) in place of the constant term
and interpreting $(x;q)_a$ etc.~according to (2.7). This gives
$$
\prod_{l=1}^n \int_{-1/2}^{1/2} dy_l \,
P_\kappa(w_1, \dots, w_n;q,q^\lambda) \prod_{i = 1}^n
(w_i; q)_a \Big ({q \over w_i};q \Big )_b
\prod_{1 \le i < j \le n} \Big ({w_i \over w_j};q \Big )_\lambda
\Big ({w_j \over w_i};q \Big )_\lambda
$$
$$
= {n! \over [n]_{q^\lambda}!} \times (\mbox{final equality in (A1)}).
$$
We now choose $a=0$, $b=\epsilon$ and apply Proposition A1 with
$$
u_\epsilon (w_l) = \big ( {q \over w_l}; q \big )_\epsilon
$$
and
$$
f(w_1, \dots, w_n) = P_\kappa(w_1, \dots, w_n;q,q^\lambda)
\prod_{1 \le i < j \le n} \Big ({w_i \over w_j};q \Big )_\lambda
\Big ({w_j \over w_i};q \Big )_\lambda
\eqno ({\rm A}3)
$$
($f$ is homogeneous of order $|\kappa|$). This gives
\renewcommand{\theequation}{A4}
\begin{eqnarray}\lefteqn{
 \prod_{l=2}^n \int_{-1/2}^{1/2}
dy_l\,f(1,w_2, \dots, w_n)} \nonumber \\ & & 
={[|\kappa|]_q \Gamma_q(\kappa_1) \over n }
{n! \over [n]_{q^\lambda}! } {\Gamma_q(\lambda n + 1) \over (\Gamma_q( 
\lambda + 1)^n}P_\kappa
(1,q^\lambda,
\dots q^{(n-1)\lambda}) {{[0]'}^{(1/\lambda)}_{\kappa;q} \over
[1 + (n-1) \lambda]_{\kappa;q}^{(1/\lambda)}},
\end{eqnarray}
where we have used the formulas
$$
\lim_{\epsilon \rightarrow 0}{1 \over
\epsilon}[-\epsilon]_{\kappa;q}^{(1/\lambda)} =  \log q { \Gamma_q(\kappa_1)
\over 1 - q}
{[0]'}^{(1/\lambda)}_{\kappa;q}
$$
and
$$
\int_{-1/2}^{1/2} dy_l \, \Big ({q \over w_l},q \Big )_\epsilon
w_l^{|\kappa|} \: \sim
\:  \epsilon \log q { q^{|\kappa|} \over 1 - q^{|\kappa|}}
$$
and we have used the fact that $D_0(0,0,\lambda;q)$ is
given by (2.5) with
$N_0 = n$.

We remark that the formula (A4) is of direct relevance to the calculation of
correlation functions in the so-called relativistic Calogero-Sutherland model
\cite{rusc}
(the calculation of these correlations has been announced by Konno \cite{konn}).
It provides the Fourier coefficients in the expansion of a symmetric sum of Dirac
delta functions in terms of Macdonald polynomials. To see this we note that
$$
\{ P_\kappa(w_1, \dots, w_n;q;q^\lambda)\}_\kappa \: \bigcup \:
\{ \prod_{j=1}^n w_j^{-l} \,P_\kappa(w_1, \dots, w_n;q;q^\lambda)
\}_{\kappa : \kappa_n = 0, \, l = 1,2,\dots} 
$$
 form a complete set of functions which are orthogonal with respect
to the inner product
$$
\langle f | g \rangle' :={1 \over n!} \prod_{l=1}^n \int_{-1/2}^{1/2} dy_l \,
\prod_{j \ne k}^n \big ( {w_j \over w_k}; q \big )_\lambda f^* g.
$$
 The Fourier formula then gives
$$
\sum_{j=1}^n \delta(y_j) = \sum_\kappa {\beta_\kappa \over
\langle P_\kappa | P_\kappa \rangle'} P_\kappa (w_1, \dots, w_n;q,q^\lambda)
+ \sum_{l=1}^\infty \sum_{\kappa: \,\kappa_n = 0} {\gamma_{l,\kappa} \over
\langle P_\kappa | P_\kappa \rangle'} \prod_{j=1}^n w_j^{-l} \,
 P_\kappa (w_1, \dots, w_n;q,q^\lambda)
\eqno ({\rm A}5)
$$
where
$$
\beta_\kappa = {n \over n!} \prod_{l=2}^n \int_{-1/2}^{1/2}
dy_l\,f(1,w_2, \dots, w_n)
$$
and similarly for $\gamma_{l,\kappa}$. Thus (A4) immediately gives the value of
$\beta_\kappa$.

To derive the formula for the power sum expansion from (A5) we note that
$$
\sum_{j=1}^n \delta(y_j) = \sum_{k= -\infty}^\infty (w_1^k + \cdots + w_n^k ).
$$
Since $ P_\kappa (w_1, \dots, w_n;q,q^\lambda)$ is homogeneous of order
$|\kappa|$,
it follows by equating terms homogenoeus of order $|\kappa|$ on both sides that
$$
w_1^k + \cdots + w_n^k = \sum_{|\kappa| = k} {\beta_\kappa \over
\langle P_\kappa | P_\kappa \rangle'} P_\kappa (w_1, \dots, w_n;q,q^\lambda).
$$
Substituting $n$ times the r.h.s.~of (A4) for $\beta_\kappa$ gives the stated
result.

\vspace{2mm}
We should remark here, that the coefficients relating the
power sums $\tilde{p}_k(w_1,\ldots,w_n):=\sum_{i=1}^n w_i^n$ and 
the Macdonald polynomials
$P_{\ka}(w_1,\ldots,w_n;q,q^{\la})$ can also be deduced from certain
results in Macdonald's book\cite{macd}. Let $t=q^{\la}$ as before, 
and for a partition $\si=(n^{f_n}\ldots 2^{f_2}1^{f_1})$, let
$$
z_{\si}(t) = \prod_{i}i^{f_i}f_i!\:(1-t^i)^{-f_i}
$$
and also define
\renewcommand{\theequation}{A6}
\begin{eqnarray}
c_{\ka}(q,t) &=& \prod_{s\in\ka}(1-q^{a(s)}t^{l(s)+1}) \nonumber\\
c_{\ka}'(q,t) &=& \prod_{s\in\ka}(1-q^{a(s)+1}t^{l(s)})
\end{eqnarray}
where $a(s)$ (respectively $l(s)$) are the the number of squares to the
right (resp. underneath) the node $s$ in the Ferrer's diagram of
$\ka$. Macdonald introduces functions $X^{\ka}_{\si}(q,t)$ 
(which are conjectured to be polynomials in $q$ and $t$) such that
for an arbitrary number of indeterminates $w_i$,
$$
P_{\ka}(w;q,t) = \frac{1}{c_{\ka}(q,t)}\sum_{\si} \frac{1}{z_{\si}(t)}
\,X^{\ka}_{\si}(q,t)\:\tilde{p}_{\si}(w)
$$
The polynomials $X^{\ka}_{\si}(q,t)$ obey an orthogonality
relation which allows us to invert the above equation, yielding
$$
\tilde{p}_{\si}(w) = \prod_i(1-q^i)^{f_i}\,\sum_{\ka} \frac{1}{c_{\ka}'(q,t)}
\,X^{\ka}_{\si}(q,t)\,P_{\ka}(w;q,t) 
$$
In the particular case $\si=(k)$, there is the explicit 
formula \cite[p 366]{macd}
$$
X^{\ka}_{(k)}(q,t) = \prod_{(i,j)\in\ka} (t^{i-1}-q^{j-1})
\eqno ({\rm A}7)
$$
where the product is over all nodes $(i,j)$ in $\ka$ (labelled
in matrix-fashion) {\it excluding} the node $(1,1)$. Macdonald
has essentially shown that (A6) and (A7)
can be re-written in ``label-dependent'' forms
\begin{eqnarray*}
X^{\ka}_{(k)}(q,t) &=& t^{n(\ka)}\,(q;q)_{\ka_1-1}\,
\prod_{i=2}^r (t^{-i+1};q)_{\ka_i} \\
c_{\ka}'(q,t) &=& \prod_{i=1}^r \frac{1}{(q;q)_{\ka_i +\la(r-i)}}
\: \prod_{1\leq i<j\leq r} (q^{\ka_i-\ka_j+1+\la(j-i-1)};q)_{\la}
\end{eqnarray*}
where $r$ is the length of $\ka$ (that is, the number of
non-zero parts), and $n(\ka):=\sum_i (i-1)\ka_i$.
It thus follows that $\tilde{p}_k(w)=\sum_{|\ka|=k}
a_{\ka}\,P_{\ka}(w;q,q^{\la})$ where
\renewcommand{\theequation}{A8}
\begin{eqnarray}
a_{\ka} = q^{\la\,n(\ka)}\,[k]_q\,\Ga_q(\ka_1)\prod_{i=2}^r
\frac{\Ga_q(\la(-i+1)+\ka_i)}{\Ga_q(\la(-i+1))} \,\prod_{i=1}^r
\frac{1}{\Ga_q(\la(r-i)+\ka_i)} \nonumber\\
\times\prod_{1\leq i<j\leq r}
\frac{\Ga_q(\ka_i-\ka_j+1+\la(j-i))}{\Ga_q(\ka_i-\ka_j+1+\la(j-i-1))}
\end{eqnarray}
Through simplification one can show that
$$
\frac{\alpha_{\ka;q}}{\langle P_{\ka} | P_{\ka}
\rangle '} = a_{\ka}
$$
thus providing an alternative proof of Proposition A2. Furthermore,
the formula (A8) explicitly demonstrates that the coefficients in 
the expansion are independent of the number of variables $n$.

\pagebreak

\end{document}